\documentclass[10pt,a4paper]{amsart}
\usepackage{amssymb,amsmath}
\usepackage[american,french]{babel}
\usepackage{amsopn}
\usepackage{graphicx}
\usepackage{fancyhdr}
\usepackage[T1]{fontenc}
\title{SUBELLIPTIC ESTIMATES FOR SYSTEMS OF QUADRATIC DIFFERENTIAL OPERATORS}
\newcommand{\rr}{\mathbb{R}}
\newcommand{\eps}{\varepsilon}
\newcommand{\nn}{\mathbb{N}}
\newcommand{\cc}{\mathbb{C}}

\def\init{\setcounter{equa}{0}}
\def\inc{\stepcounter{equa}}
\def\num{\tag{\thesection.\theequa}}

\def\wrtext#1{\relax\ifmmode{\leavevmode\hbox{#1}}\else{#1}\fi}

\def\begeq{\begin{equation}}
\def\endeq{\end{equation}}

\def\part#1{\frac{\partial}{\partial #1}}

\def\og{orthogonal}

\begin{document}
\newcounter{equa}
\selectlanguage{american}

\newtheorem{lemma}{Lemma}[subsection]
\newtheorem{definition}{Definition}[subsection]
\newtheorem{proposition}{Proposition}[subsection]
\newtheorem{theorem}{Theorem}[subsection]

\title[SUBELLIPTIC ESTIMATES FOR SYSTEMS OF QUADRATIC OPERATORS]{SUBELLIPTIC ESTIMATES FOR OVERDETERMINED SYSTEMS OF QUADRATIC DIFFERENTIAL OPERATORS}
\author{Karel \textsc{Pravda-Starov}}
\address{\noindent \textsc{Department of Mathematics,
Imperial College London,
Huxley Building, 180 Queen's Gate, 
London SW7 2AZ, UK}
} 
\email{k.pravda-starov@imperial.ac.uk}
\urladdr{http://www2.imperial.ac.uk/~kpravdas/index.html}

\begin{abstract}
We prove global subelliptic estimates for systems of quadratic differential operators. Quadratic differential operators are operators defined in the Weyl quantization by complex-valued quadratic symbols. In a previous work, we pointed out the existence of a particular linear subvector space in the phase space intrinsically associated to their Weyl symbols, 
called singular space, which rules a number of fairly general properties of non-elliptic quadratic operators. About the subelliptic properties of these operators, we established that quadratic operators with zero singular spaces fulfill global subelliptic estimates with a loss of derivatives depending on certain algebraic properties of the Hamilton maps associated to their Weyl symbols. 
The purpose of the present work is to prove similar global subelliptic estimates for overdetermined systems of quadratic operators. We establish here a simple criterion for the subellipticity of these systems giving an explicit measure of the loss of derivatives and highlighting the non-trivial interactions played by the different operators composing those systems. 
\end{abstract}


\maketitle

\section{Introduction}
\init

\subsection{Miscellaneous facts about quadratic differential operators}
In a recent joint work with M.~Hitrik, we investigated spectral and semigroup properties of non-elliptic quadratic operators. 
Quadratic operators are pseudodifferential operators
defined in the Weyl quantization
\begin{equation}\label{3}\inc
q^w(x,D_x) u(x) =\frac{1}{(2\pi)^n}\int_{\rr^{2n}}{e^{i(x-y).\xi}q\Big(\frac{x+y}{2},\xi\Big)u(y)dyd\xi}, \num
\end{equation}
by some symbols $q(x,\xi)$, with $(x,\xi) \in \rr^{n} \times \rr^n$ and $n \in \nn^*$, which
are complex-valued quadratic forms. Since these symbols are quadratic forms, the corresponding operators in
(\ref{3}) are in fact differential operators. Indeed, the Weyl quantization of the quadratic symbol
$x^{\alpha} \xi^{\beta}$, with $(\alpha,\beta) \in \nn^{2n}$ and $|\alpha+\beta| = 2$, is the differential operator
$$
\frac{x^{\alpha}D_x^{\beta}+D_x^{\beta} x^{\alpha}}{2}, \ D_x=i^{-1}\partial_x.
$$
One can also notice that quadratic differential operators
are a priori formally non-selfadjoint since their Weyl symbols in (\ref{3}) are complex-valued.

\medskip

Considering quadratic operators whose Weyl symbols have real parts with a sign, say here, Weyl symbols with non-negative real parts 
\begin{equation}\label{smm1}\inc
\textrm{Re }q \geq 0, \num
\end{equation}
we pointed out in \cite{hps} the existence of a particular linear subvector space $S$ in the phase space $\rr_x^n \times \rr_{\xi}^n$ intrinsically associated to their Weyl symbols $q(x,\xi)$, 
called singular space, which seems to play a basic r\^ole in the understanding of a number of fairly general properties of non-elliptic quadratic operators.
More specifically, we first proved in \cite{hps} (Theorem~1.2.1) that when the singular space $S$
has a symplectic structure then the associated heat equation
\begin{equation}\label{sm1}\inc
\left\lbrace
\begin{array}{c}
\displaystyle \frac{\partial u}{\partial t}(t,x)+q^w(x,D_x) u(t,x)=0  \\
u(t,\textrm{\textperiodcentered})|_{t=0}=u_0 \in L^2(\rr^n),
\end{array} \right.
 \num
\end{equation}
is smoothing in every direction of the orthogonal complement $S^{\sigma \perp}$ of $S$ with respect to the canonical
symplectic form $\sigma$ on $\rr^{2n}$,
\begin{equation}
\label{11}\inc
\sigma \big{(}(x,\xi),(y,\eta) \big{)}=\xi.y-x.\eta, \ (x,\xi) \in \rr^{2n},  (y,\eta) \in \rr^{2n}, \num
\end{equation}
that is, that, if $(x',\xi')$ are some linear symplectic coordinates on the symplectic space $S^{\sigma \perp}$
then we have for all $t>0$, $N \in \nn$ and $u \in L^2(\rr^n)$,
\begin{equation}\label{sm4bis}\inc
\big((1+|x'|^2+|\xi'|^2)^N\big)^we^{-t q^w(x,D_x)}u \in L^2(\rr^n). \num
\end{equation}
We also proved in \cite{hps} (See Section~1.4.1 and Theorem~1.2.2) that when the Weyl symbol $q$ of a quadratic operator fulfills (\ref{smm1}) and an assumption of
partial ellipticity on its singular space $S$ in the sense that
\begin{equation}\label{sm2}\inc
(x,\xi) \in S, \ q(x,\xi)=0 \Rightarrow (x,\xi)=0, \num
\end{equation}
then this singular space always has a symplectic structure and the spectrum of the operator
$q^w(x,D_x)$ is only composed of a countable number of eigenvalues of finite multiplicity, with a similar structure
as the one established by J. Sj\"ostrand for elliptic quadratic operators in his classical work~\cite{sjostrand}. Elliptic quadratic operators are the quadratic operators whose symbols satisfy the condition of global ellipticity
$$(x,\xi) \in \rr^{2n}, \ q(x,\xi)=0 \Rightarrow (x,\xi)=0,$$
on the whole phase space $\rr^{2n}$. Let us recall here that spectral properties of quadratic operators are playing a basic r\^ole in the analysis of partial differential operators with double characteristics. This is particularly the case in some general results about hypoellipticity. We refer the reader to~\cite{hypoelliptic}, \cite{sjostrand}, as well as Chapter 22 of~\cite{hormander} together with all the references given there.

\medskip

In the present paper, we are interested in studying the subelliptic properties of overdetermined systems of non-selfadjoint quadratic operators. 
This work can be viewed as a natural extension of the analysis led in \cite{sub}, in which we investigated  in the scalar case the r\^ole played by the singular space when studying subelliptic properties of quadratic operators. We aim here at showing how the analysis led in this previous work can be pushed further when dealing with overdetermined systems of quadratic operators. We shall see that the techniques introduced in \cite{sub} are sufficiently robust to be extended to the system case and that they turn out to be sufficiently sharp to highlight phenomena of non-trivial interactions between the different quadratic operators composing a system. In this paper, we shall therefore be interested in establishing some global subelliptic estimates of the type
\begin{equation}\label{lol1}\inc
\big\|\big(\langle(x,\xi)\rangle^{2(1-\delta)}\big)^w u\big\|_{L^2} \lesssim \sum_{j=1}^N\|q_j^w(x,D_x) u\|_{L^2}+\|u\|_{L^2},\num
\end{equation}
where $\langle(x,\xi)\rangle=(1+|x|^2+|\xi|^2)^{1/2}$ and $\delta>0$; for systems of the $N$ quadratic operators 
$q_j^w(x,D_x)$, with $1 \leq j \leq N$. The positive parameter $\delta>0$ appearing in (\ref{lol1}) will measure the loss of derivatives  with respect to the elliptic case (case $\delta=0$). As in the scalar case studied in~\cite{sub}, we aim at giving a simple criterion for systems of quadratic operators ensuring that a global subelliptic estimate of the type (\ref{lol1}) holds together with an explicit characterization of the associated loss of derivatives. This loss of derivatives $\delta$ will be characterized in terms of algebraic conditions on the Hamilton maps associated to the Weyl symbols of the quadratic operators composing the system.

\medskip

In this work, we study the subellipticity of overdetermined systems in the sense given by P.~Bolley, J.~Camus and J.~Nourrigat in \cite{Camus} (Theorem~1.1). In this seminal work, these authors study the microlocal subellipticity of overdetermined systems of pseudodifferential operators. More specifically, they establish the subellipticity of systems composed of pseudodifferential operators with real principal symbols satisfying the H\"ormander-Kohn condition. More generally, in the case of overdetermined systems of non-selfadjoint pseudodifferential operators, the greatest achievements up to now were obtained by J.~Nourrigat in~\cite{Nou1} and~\cite{Nou2}. In these two major works, J.~Nourrigat studies the microlocal subellipticity and maximal hypoellipticity for systems of non-selfadjoint pseudodifferential operators by the mean of representations of nilpotent groups. We shall explain in the following how the algebraic condition on the Hamilton maps  (\ref{lid1}) in Theorem~\ref{theorem1} relates with these former results. More specifically, we shall comment on its link with the H\"ormander-Kohn condition appearing in~\cite{Camus} (Theorem~1.1).

\medskip

Before giving the precise statement of our main result, we shall recall miscellaneous notations about quadratic differential operators and the results obtained in the scalar case. In all the following, we consider
\begin{eqnarray*}
q_j : \rr_x^n \times \rr_{\xi}^n &\rightarrow& \cc\\
 (x,\xi) & \mapsto & q_j(x,\xi),
\end{eqnarray*}
with $1 \leq j \leq N$, $N$ complex-valued quadratic forms with non-negative real parts
\begin{equation}\label{inf1}\inc
\textrm{Re }q_j(x,\xi) \geq 0, \ (x,\xi) \in \rr^{2n}, n \in \nn^*. \num
\end{equation}
We know from \cite{mehler} (p.425) that the
maximal closed realization of a quadratic operator $q^w(x,D_x)$ whose Weyl symbol has a non-negative real part, i.e., the operator on $L^2(\rr^n)$ with the domain
$$D(q)=\big\{u \in L^2(\rr^n) :  q^w(x,D_x) u \in L^2(\rr^n)\big\},$$
coincides with the graph closure of its restriction to $\mathcal{S}(\rr^n)$,
$$q^w(x,D_x) : \mathcal{S}(\rr^n) \rightarrow \mathcal{S}(\rr^n).$$
Associated to a quadratic symbol $q$ is the numerical range $\Sigma(q)$ defined as the closure in the
complex plane of all its values
\begin{equation}\label{9}\inc
\Sigma(q)=\overline{q(\rr_x^n \times \rr_{\xi}^n)}. \num
\end{equation}
We also recall from~\cite{hormander} that the Hamilton map $F \in M_{2n}(\cc)$ associated to the quadratic form $q$
is the map uniquely defined by the identity
\begin{equation}\label{10}\inc
q\big{(}(x,\xi);(y,\eta) \big{)}=\sigma \big{(}(x,\xi),F(y,\eta) \big{)}, \ (x,\xi) \in \rr^{2n},  (y,\eta) \in \rr^{2n}, \num
\end{equation}
where $q\big{(}\textrm{\textperiodcentered};\textrm{\textperiodcentered} \big{)}$ stands for the polarized form
associated to the quadratic form $q$. 
It directly follows from the definition of the Hamilton map $F$ that
its real part and its imaginary part
$$\textrm{Re } F=\frac{1}{2}(F+\overline{F}) \textrm{ and } \textrm{Im }F=\frac{1}{2i}(F-\overline{F}),$$ 
are the Hamilton maps associated
to the quadratic forms $\textrm{Re } q$ and $\textrm{Im }q$, respectively. One can also notice from (\ref{10}) that an
Hamilton map is always skew-symmetric with respect to $\sigma$. This is just a consequence of the
properties of skew-symmetry of the symplectic form and symmetry of the polarized form
\begin{equation}\label{12}\inc
\forall X,Y \in \rr^{2n}, \ \sigma(X,FY)=q(X;Y)=q(Y;X)=\sigma(Y,FX)=-\sigma(FX,Y).\num
\end{equation}

Associated to the symbol $q$, we defined in~\cite{hps} its singular space $S$ as
the following intersection of kernels
\begin{equation}\label{h1}\inc
S=\Big(\bigcap_{j=0}^{+\infty}\textrm{Ker}\big[\textrm{Re }F(\textrm{Im }F)^j \big]\Big) \cap \rr^{2n}, \num
\end{equation}
where the notations $\textrm{Re } F$ and $\textrm{Im }F$ stand respectively for the real part and the imaginary part of the Hamilton map associated to $q$.
Notice that the Cayley-Hamilton theorem applied to $\textrm{Im }F$ shows that
$$(\textrm{Im }F)^k X \in \textrm{Vect}\big(X,...,(\textrm{Im }F)^{2n-1}X\big), \ X \in \rr^{2n}, \ k \in \nn,$$
where $\textrm{Vect}\big(X,...,(\textrm{Im }F)^{2n-1}X\big)$ is the vector space spanned by the vectors $X$, ...,
$(\textrm{Im }F)^{2n-1}X$; and therefore
the singular space is actually equal to the following finite intersection of the kernels
\begin{equation}\label{h1bis}\inc
S=\Big(\bigcap_{j=0}^{2n-1}\textrm{Ker}\big[\textrm{Re }F(\textrm{Im }F)^j \big]\Big) \cap \rr^{2n}. \num
\end{equation}

Considering a quadratic operator $q^w(x,D_x)$ whose Weyl symbol  
\begin{eqnarray*}
q : \rr_x^n \times \rr_{\xi}^n &\rightarrow& \cc\\
 (x,\xi) & \mapsto & q(x,\xi),
\end{eqnarray*}
has a non-negative real part, $\textrm{Re }q \geq 0$, we established in \cite{sub} (Theorem~1.2.1) that when its singular space $S$ is reduced to $\{0\}$, the operator $q^w(x,D_x)$ fulfills the following global subelliptic estimate 
\begin{equation}\label{kristen1}\inc
\exists C>0, \forall u \in D(q), \  \big\|\big(\langle(x,\xi)\rangle^{2/(2k_0+1)}\big)^w u\big\|_{L^2} \leq C\big(\|q^w(x,D_x) u\|_{L^2}+\|u\|_{L^2}\big), \num
\end{equation}
where $k_0$ stands for the smallest non-negative integer, $0 \leq k_0 \leq 2n-1$, such that the intersection of the following $k_0+1$ kernels with the phase space $\rr^{2n}$ is reduced to~$\{0\}$,
\begin{equation}\label{h1biskrist2}\inc
\Big(\bigcap_{j=0}^{k_0}\textrm{Ker}\big[\textrm{Re }F(\textrm{Im }F)^j \big]\Big) \cap \rr^{2n}=\{0\}. \num
\end{equation}
Notice that the loss of derivatives $\delta=2k_0/(2k_0+1)$, appearing in the subelliptic estimate (\ref{kristen1}) directly depends on the non-negative integer $k_0$ characterized by the algebraic condition (\ref{h1biskrist2}).

More generally, considering a quadratic operator $q^w(x,D_x)$ whose Weyl symbol has a non-negative real part with a singular space $S$ which may differ from $\{0\}$, but does have a symplectic structure in the sense that the restriction of the canonical symplectic form $\sigma$ to $S$ is non-degenerate, we proved in \cite{sub} (Theorem~1.2.2) that 
the operator $q^w(x,D_x)$ is subelliptic in any direction of the orthogonal complement $S^{\sigma \perp}$ of the singular space with respect to the symplectic form $\sigma$ in the sense that, if $(x',\xi')$ are some linear symplectic coordinates on $S^{\sigma \perp}$ then we have 
$$\exists C>0, \forall u \in D(q), \  \big\|\big(\langle(x',\xi')\rangle^{2/(2k_0+1)}\big)^w u\big\|_{L^2} \leq C\big(\|q^w(x,D_x) u\|_{L^2}+\|u\|_{L^2}\big),$$
with $\langle(x',\xi')\rangle=(1+|x'|^2+|\xi'|^2)^{1/2}$, where $k_0$ stands for the smallest non-negative integer, $0 \leq k_0 \leq 2n-1$, such that \
\begin{equation}\label{h1biskrist2bis}\inc
S=\Big(\bigcap_{j=0}^{k_0}\textrm{Ker}\big[\textrm{Re }F(\textrm{Im }F)^j \big]\Big) \cap \rr^{2n}. \num
\end{equation}
Finally, we end these few recalls by underlining that the assumption about the symplectic structure of the singular space is always fulfilled by any quadratic symbol $q$ which satisfies the assumption of partial ellipticity on its singular space $S$,
$$(x,\xi) \in S, \ q(x,\xi)=0 \Rightarrow (x,\xi)=0.$$
We refer the reader to Section~1.4.1 in~\cite{hps} for a proof of this fact.

\subsection{Statement of the main result}

Considering a system of $N$ quadratic operators $q_j^w(x,D_x)$, $1 \leq j \leq N$, whose Weyl symbols $q_j$ have all non-negative real parts
\begin{equation}\label{inf1imp}\inc
\textrm{Re }q_j(x,\xi) \geq 0, \ (x,\xi) \in \rr^{2n}, \ n \in \nn^*, \num
\end{equation}
and denoting by $F_j$ their associated Hamilton maps, the main result contained in this article is the following: 

\bigskip

\begin{theorem}\label{theorem1}
Consider a system of $N$ quadratic operators $q_j^w(x,D_x)$, $1 \leq j \leq N$, satisfying \emph{(\ref{inf1imp})}. If there exists $k_0 \in \nn$ such that
\begin{equation}\label{lid1}\inc 
\Big(\bigcap_{0 \leq k \leq k_0} \bigcap_{\substack{j=1,...,N, \\ 
(l_1,...,l_{k}) \in \{1,...,N\}^k}
}\emph{\textrm{Ker}}(\emph{\textrm{Re }}F_j\emph{\textrm{Im }}F_{l_1}...\emph{\textrm{Im }}F_{l_{k}})\Big) \cap \rr^{2n}=\{0\},\num
\end{equation}
then this overdetermined system of quadratic operators is subelliptic with a loss of $\delta=2k_0/(2k_0+1)$ derivatives, that is, that there exists $C>0$ such that for all $u \in D(q_1) \cap ... \cap D(q_N)$,
\begin{equation}\label{kristen5}\inc
\big\|\big(\langle(x,\xi)\rangle^{2/(2k_0+1)}\big)^w u\big\|_{L^2} \leq C\Big(\sum_{j=1}^N\|q_j^w(x,D_x) u\|_{L^2}+\|u\|_{L^2}\Big), \num
\end{equation}
with $\langle(x,\xi)\rangle=(1+|x|^2+|\xi|^2)^{1/2}$.
\end{theorem}

\bigskip

\noindent \textit{Remark.} Let us make clear that the intersection of kernels  
$$\bigcap_{\substack{j=1,...,N, \\  (l_1,...,l_{k}) \in \{1,...,N\}^k}}\textrm{Ker}(\textrm{Re }F_j\textrm{Im }F_{l_1}...\textrm{Im }F_{l_{k}}),$$
is to be understood as 
$$\bigcap_{j=1,...,N}\textrm{Ker} \ \textrm{Re }F_j,$$
when $k=0$.

\subsection{Example of a subelliptic system of quadratic operators} The following example of subelliptic system of quadratic operators shows that Theorem~\ref{theorem1} really highlights new non-trivial interaction phenomena between the different operators composing a system, which cannot be derived from the result of subellipticity known in the scalar case  (Theorem~1.2.1 in \cite{sub}). Indeed, define the quadratic forms
$$q_j(x,\xi)=x_1^2+\xi_1^2+i(\xi_1^2+x_{j+1} \xi_1) \textrm{ and } \tilde{q}_j(x,\xi)=x_1^2+\xi_1^2+i(\xi_1^2+\xi_{j+1} \xi_1),$$
for $1 \leq j \leq n-1$ and $(x,\xi) \in \rr^{2n}$, with $n \geq 2$. A direct computation using (\ref{10}) and (\ref{h1bis}) shows that the singular space of the quadratic form
$$\sum_{j=1}^{n-1}(\lambda_j q_j+\tilde{\lambda}_j\tilde{q}_j),$$
for some real numbers $\lambda_j,\tilde{\lambda}_j$ verifying 
$$\sum_{j=1}^{n-1}(\lambda_j+\tilde{\lambda}_j)> 0;$$ is given by 
$$S=\Big\{(x,\xi) \in \rr^{2n} : x_1=\xi_1=\sum_{j=1}^{n-1}(\lambda_j x_{j+1}+\tilde{\lambda}_j\xi_{j+1})=0\Big\},$$
which is always a non-zero subvector space. It then follows that one cannot deduce any result about the subellipticity of the scalar operator
$$\sum_{j=1}^{n-1}(\lambda_j q_j^w(x,D_x)+\tilde{\lambda}_j\tilde{q}_j^w(x,D_x)),$$
in order to get the subellipticity of the overdetermined system composed by the $2n-2$ operators $q_j^w(x,D_x)$ and $\tilde{q}_j^w(x,D_x)$, for  $1 \leq j \leq n-1$. Nevertheless, by denoting respectively $F_j$ and $\tilde{F}_j$ the Hamilton maps of the quadratic forms $q_j$ and $\tilde{q}_j$, another direct computation using (\ref{10}) shows that 
$$\textrm{Ker }\textrm{Re }F_j \cap \textrm{Ker}(\textrm{Re }F_j\textrm{Im }F_{j}) \cap \rr^{2n}=\{(x,\xi) \in \rr^{2n} : x_1=\xi_1=x_{j+1}=0\}$$
and
$$\textrm{Ker }\textrm{Re }\tilde{F}_j \cap \textrm{Ker}(\textrm{Re }\tilde{F}_j\textrm{Im }\tilde{F}_{j}) \cap \rr^{2n}=\{(x,\xi) \in \rr^{2n} : x_1=\xi_1=\xi_{j+1}=0\}.$$
One can then deduce from Theorem~\ref{theorem1} the following global subelliptic estimate with a loss of $2/3$ derivatives 
$$\big\|\big(\langle(x,\xi)\rangle^{2/3}\big)^w u\big\|_{L^2} \lesssim \sum_{j=1}^{n-1}\big(\|q_j^w(x,D_x) u\|_{L^2}+\|\tilde{q}_j^w(x,D_x) u\|_{L^2}\big) +\|u\|_{L^2}.$$
Of course, Theorem~\ref{theorem1} can highlight more complex interactions between the different operators composing the system when we consider operators with different real parts.

 \subsection{Comments on the condition for subellipticity}
Theorem~\ref{theorem1} gives a very explicit and simple algebraic condition on the Hamilton maps of quadratic operators ensuring the subellipticity of the system. Let us notice that this condition is very easy to handle and allows to directly measure the associated loss of derivatives by a straightforward computation. 
We shall now explain how this is related to the H\"ormander-Kohn condition. Recall from~\cite{Camus} (Theorem~1.1) that the H\"ormander-Kohn condition for microlocal subellipticity of overdetermined systems of pseudodifferential operators with real principal symbols; reads as the existence of an elliptic iterated commutator of the operators composing the system. In the case of a system of non-selfadjoint quadratic operators $(q_j^w)_{1 \leq j \leq N}$, if we assume in addition that this system is maximal hypoelliptic\footnote{We refer to \cite{Nou1} and \cite{Nou2} for conditions and general results of maximal hypoellipticity for overdetermined systems of non-selfadjoint pseudodifferential operators.}, the natural condition becomes to ask the ellipticity of an iterated commutator of the real parts $((\textrm{Re } q_j)^w)_{1 \leq j \leq N}$ and imaginary parts $((\textrm{Im } q_j)^w)_{1 \leq j \leq N}$ of the operators composing the system.  
Coming back to our specific condition for subellipticity (\ref{lid1}), we first notice that in the scalar case, it reads as the existence of a non-negative integer $k_0$ such that
$$\Big(\bigcap_{j=0}^{k_0}\textrm{Ker}[\textrm{Re }F(\textrm{Im }F)^j]\Big) \cap \rr^{2n}=\{0\},$$
with $F$ standing for the Hamilton map of the unique operator $q^w(x,D_x)$ composing the system. As recalled in~\cite{sub} (Section~1.2), this condition implies that, for any non-zero point in the phase space $X_0 \in \rr^{2n}$, we can find a non-negative integer $k$ such that 
$$\forall \ 0 \leq j \leq 2k-1, \ H_{\textrm{Im}q}^j \textrm{Re }q(X_0)=0 \textrm{ and } H_{\textrm{Im}q}^{2k} \textrm{Re }q(X_0) \neq 0,$$
where $H_{\textrm{Im}q}$ stands for the Hamilton vector field of $\textrm{Im }q$,
$$H_{\textrm{Im}q}=\frac{\partial \textrm{Im }q}{\partial \xi}\cdot \frac{\partial}{\partial x}-\frac{\partial \textrm{Im }q}{\partial x} \cdot \frac{\partial}{\partial \xi}.$$
This shows that the $2k^{\textrm{th}}$ iterated commutator 
$$[\textrm{Im }q^w,[\textrm{Im }q^w,[...,[\textrm{Im }q^w,\textrm{Re }q^w]]]...]=(-1)^k(H_{\textrm{Im}q}^{2k} \textrm{Re }q)^w,$$
with exactly $2k$ terms $\textrm{Im }q^w$ in left-hand-side of the above formula; is elliptic at $X_0$; and underlines the intimate link between (\ref{lid1}) and the H\"ormander-Kohn condition in the scalar case. In the system case, the situation is more complicated and this link is less obvious to highlight explicitly. More specifically, we shall see in this case that the algebraic condition (\ref{lid1}) implies that the quadratic form
$$\sum_{k=0}^{k_0}\sum_{\substack{j=1,...,N, \\
(l_1,...,l_{k}) \in \{1,...,N\}^k}}{\textrm{Re }q_j(\textrm{Im }F_{l_1}... \textrm{Im }F_{l_{k}}X)},$$
is positive definite. This property implies that for any non-zero point $X_0 \in \rr^{2n}$, one can find $k \in \nn$, $j \in \{1,...,N\}$ and $(l_1,...,l_{k}) \in \{1,...,N\}^k$ such that 
$$\textrm{Re }q_j(\textrm{Im }F_{l_1}... \textrm{Im }F_{l_{k}}X_0)>0.$$
By considering the minimal non-negative integer $k$ with this property and using the same arguments as the ones developed in \cite{hps} (p.820-822), one can actually check that any iterated commutator of order less or equal to $2k-1$, that is,
$$[P_1,[P_2,[P_3,[...,[P_{r},P_{r+1}]...]]]],$$
with $r \leq 2k-1$, $P_l=\textrm{Re }q_{s_1}^w$ or $P_l=\textrm{Im }q_{s_2}^w$; and where at least one $P_{l_0}$ is equal to $\textrm{Re }q_{s_3}^w$, for $1\leq s_1,s_2,s_3 \leq N$; are not elliptic at $X_0$. One can also check that the non-zero term 
$$\textrm{Re }q_j(\textrm{Im }F_{l_1}... \textrm{Im }F_{l_{k}}X_0)>0,$$ 
actually appears when expanding the Weyl symbol at $X_0$ of the $2k^{\textrm{th}}$ iterated commutator
\begin{multline*}
[\textrm{Im }q_{l_k}^w,[\textrm{Im }q_{l_k}^w,[\textrm{Im }q_{l_{k-1}}^w,[\textrm{Im }q_{l_{k-1}}^w,[...,[\textrm{Im }q_{l_1}^w,[\textrm{Im }q_{l_1}^w,\textrm{Re }q_j^w]]]...]\\
=(-1)^k(H_{\textrm{Im}q_{l_k}}^{2}...H_{\textrm{Im}q_{l_1}}^{2} \textrm{Re }q_j)^w.
\end{multline*}
However, contrary to the scalar case, there may be also other non-zero terms in this expansion; and it is not really clear if this natural commutator associated to the term 
$$\textrm{Re }q_j(\textrm{Im }F_{l_1}... \textrm{Im }F_{l_{k}}X_0),$$
is actually elliptic at $X_0$,
$$H_{\textrm{Im}q_{l_k}}^{2}...H_{\textrm{Im}q_{l_1}}^{2} \textrm{Re }q_j(X_0) \overset{?}{\neq} 0.$$
Though it may be difficult to determine exactly at each point which specific commutator is elliptic, it is very likely that condition (\ref{lid1}) ensures that the H\"ormander-Kohn condition is fulfilled at any non-zero point of the phase space; and that these associated elliptic commutators are all of order less or equal to $2k_0$. It is actually what the loss of derivatives appearing in the estimate (\ref{kristen5}) suggests; and this in agreement with the optimal loss of derivatives obtained in~\cite{Camus} (Theorem~1.1) for $2k_0$ commutators 
$$\delta=1-\frac{1}{2k_0+1}=\frac{2k_0}{2k_0+1};$$
since we measure the loss of derivatives $\delta$ with respect to the elliptic case as
$$\big\|(\Lambda^{2(1-\delta)})^w u\big\|_{L^2} \lesssim \sum_{j=1}^N\|q_j^w(x,D_x) u\|_{L^2}+\|u\|_{L^2},$$
with $\Lambda^2=\langle(x,\xi)\rangle^{2}$, because quadratic operators have their Weyl symbols in the symbol class $S(\Lambda^2,\Lambda^{-2}dX^2)$ whose gain is $\Lambda^2$.

Because of the simplicity of its assumptions, Theorem~\ref{theorem1} provides a neat setting for proving global subelliptic estimates for systems of quadratic operators. It is possible that some of these global subelliptic estimates for systems of quadratic operators may also be derived from the results of microlocal subellipticity and maximal hypoellipticity proved in~\cite{Camus}, \cite{Nou1} and \cite{Nou2}. However, given a particular system of quadratic operators, one can notice that only checking the H\"ormander-Kohn condition in every non-zero point turns out to be quite difficult to do in practice. The same comment applies for checking the maximal hypoellipticity of the system. Another interest of the approach we are developing here comes from the fact that the proof of Theorem~\ref{theorem1} is purely analytic and does not require any techniques of representations of 
nilpotent groups as in~\cite{Nou1} or~\cite{Nou2}. Moreover, despite its length, the proof provided here only involves fairly elementary arguments whose complexity has no degree of comparison with the analysis led in \cite{Nou1} and \cite{Nou2}.

Finally, let us end this introduction by mentioning that this result of subellipticity for systems of quadratic operators may broaden new perspectives in the understanding of overdetermined systems of pseudodifferential operators with double characteristics; and that the construction of the weight functions in Proposition~\ref{prop1} may be of further interest and direct use in future analysis of doubly characteristic problems. In the scalar case, this construction of the weight function specific to the structure of the double characteristics obtained in~\cite{sub} (Proposition~2.0.1) has already allowed to derive in~\cite{discrete} the precise asymptotics for the resolvent norm of certain class of semiclassical pseudodifferential operators in a neighborhood of the doubly characteristic set. On the other hand, this deeper understanding of non-trivial interactions between the different quadratic operators composing overdetermined systems may also give hints on how to analyze the more complex case of $N$ by $N$ systems of quadratic operators, which is a topic of current interest. On that subject, we refer the reader to the series of recent works on non-commutative harmonic oscillators by A.~Parmeggiani and M.~Wakayama in \cite{Par6}, \cite{Par5}, \cite{Par4}, \cite{Par3}, \cite{Par2} and \cite{Par1}.

\section{Proof of Theorem~\ref{theorem1}}
\init

In the following, we shall use the notation $S_{\Omega}\big(m(X)^r,m(X)^{-2s}dX^2\big)$, where $\Omega$ is an open set in $\rr^{2n}$, $r,s \in \rr$ and $m \in C^{\infty}(\Omega,\rr_+^*)$, to stand for 
the class of symbols $a$ verifying
$$a \in C^{\infty}(\Omega), \ \forall  \alpha \in \nn^{2n}, \exists C_{\alpha}>0, \ |\partial_X^{\alpha}a(X)| \leq C_{\alpha} m(X)^{r-s|\alpha|}, \ X \in \Omega.$$
In the case where $\Omega=\rr^{2n}$, we shall drop the index $\Omega$ for simplicity.
We shall also use the notations $f \lesssim g$ and $f \sim g$, on $\Omega$, for respectively the estimates $\exists C>0$, $f \leq Cg$ and, $f \lesssim g$ and $g \lesssim f$, on $\Omega$.

The proof of Theorem~\ref{theorem1} will rely on the following key proposition. Considering for $1 \leq j \leq N$, 
\begin{eqnarray*}
q_j : \rr_x^n \times \rr_{\xi}^n &\rightarrow& \cc\\
 (x,\xi) & \mapsto & q_j(x,\xi),
\end{eqnarray*}
with $n \in \nn^*$, $N$ complex-valued quadratic forms with non-negative real parts
\begin{equation}\label{giu1}\inc
\textrm{Re }q_j(x,\xi) \geq 0, \ (x,\xi) \in \rr^{2n}, \ 1 \leq j \leq N,  \num
\end{equation}
we assume that there exist a positive integer $m \in \nn^*$ and an open set $\Omega_0$ in $\rr^{2n}$ such that the following sum of non-negative quadratic forms satisfies
\begin{equation}\label{giu2}\inc  
\exists c_0>0, \forall X \in \Omega_0, \ \sum_{k=0}^m \sum_{\substack{j=1,...,N, \\  (l_1,...,l_{k}) \in \{1,...,N\}^k}}{\textrm{Re }q_j(\textrm{Im }F_{l_1}...\textrm{Im }F_{l_k} X)} \geq c_0 |X|^2, \num
\end{equation}
where the notation $\textrm{Im }F_j$ stands for the imaginary part of the Hamilton map $F_j$ associated to the quadratic form $q_j$. Under this assumption, one can then extend the construction of the bounded weight function done in the scalar case in \cite{sub} (Proposition~2.0.1) to the system case as follows:

\bigskip

\begin{proposition}\label{prop1}
If $(q_j)_{1 \leq j \leq N}$ are $N$ complex-valued quadratic forms on $\rr^{2n}$ verifying \emph{(\ref{giu1})} and \emph{(\ref{giu2})} then there exist $N$ real-valued weight functions 
$$g_j \in S_{\Omega_0}\big(1,\langle X \rangle^{-\frac{2}{2m+1}}dX^2\big), \ 1 \leq j \leq N,$$ such that 
\begin{equation}\label{giu3}\inc
\exists c, c_1, ..., c_N>0, \forall X \in \Omega_0, \ 1+\sum_{j=1}^{N}\big(\emph{\textrm{Re }}q_j(X)+c_jH_{\emph{\textrm{Im}}q_j\ } g_j(X)\big) \geq c \langle X \rangle^{\frac{2}{2m+1}}, \num
\end{equation}
where the notation $H_{\emph{\textrm{Im}}q_j}$ stands for the Hamilton vector field of the imaginary part of~$q_j$.
\end{proposition}

\bigskip

As in \cite{sub}, the construction of these weight functions will be really the core of this work. This construction will be an adaptation to the system case of the one performed in the scalar case.

\medskip

To check that we can actually deduce Theorem~\ref{theorem1} from Proposition~\ref{prop1}, we begin by
noticing,  as in~\cite{sub}, that the assumptions of Theorem~\ref{theorem1} imply that the following sum of non-negative quadratic forms
\begin{equation}\label{giu6}\inc  
\exists c_0>0,  \ r(X)=\sum_{k=0}^{k_0}\sum_{\substack{j=1,...,N, \\
(l_1,...,l_{k}) \in \{1,...,N\}^k}}{\textrm{Re }q_j(\textrm{Im }F_{l_1}... \textrm{Im }F_{l_{k}}X)} \geq c_0 |X|^2, \num
\end{equation}
is actually a positive definite quadratic form. Let us indeed consider $X_0 \in \rr^{2n}$ such that $r(X_0)=0$. 
Then, the non-negativity of quadratic forms $\textrm{Re } q_j$ induces that for all $0 \leq k \leq k_0$, $j=1,...,N$ and $(l_1,...,l_{k}) \in \{1,...,N\}^k$,
\begin{equation}\label{giu7}\inc
\textrm{Re }q_j(\textrm{Im } F_{l_1}...\textrm{Im } F_{l_{k}} X_0)=0. \num
\end{equation}
By denoting $\textrm{Re }q_j(X;Y)$ the polar form associated to $\textrm{Re }q_j$, we deduce from the Cauchy-Schwarz inequality, (\ref{10}) and (\ref{giu7}) that
for all $Y \in \rr^{2n}$, 
\begin{align*}
|\textrm{Re }q_j(Y;\textrm{Im }F_{l_1}...\textrm{Im }F_{l_{k}} X_0)|^2= & \ |\sigma(Y,\textrm{Re }F_j \textrm{Im }F_{l_1}...\textrm{Im }F_{l_{k}} X_0)|^2 \\
 \leq & \ \textrm{Re }q_j(Y) \ \textrm{Re }q_j(\textrm{Im } F_{l_1}...\textrm{Im } F_{l_{k}} X_0) =0. 
\end{align*} 
It follows that for all $Y \in \rr^{2n}$, 
$$\sigma(Y,\textrm{Re }F_j \textrm{Im }F_{l_1}...\textrm{Im }F_{l_{k}} X_0)=0,$$
which implies that for all $0 \leq k \leq k_0$, $j=1,...,N$ and $(l_1,...,l_{k}) \in \{1,...,N\}^k$,
\begin{equation}\label{2.3.100}\inc
\textrm{Re }F_j \textrm{Im }F_{l_1}...\textrm{Im }F_{l_{k}} X_0=0, \num 
\end{equation}
since $\sigma$ is non-degenerate. We finally deduce (\ref{giu6}) from the assumption (\ref{lid1}).

In the case where $k_0=0$, we notice that the quadratic form 
$$q=q_1+...+q_N,$$ 
has a positive definite real part. This implies in particular that $q$ is elliptic on $\rr^{2n}$. One can therefore directly deduce from classical results about elliptic quadratic differential operators proved in \cite{sjostrand} (See Theorem~3.5 in \cite{sjostrand} or comments about the elliptic case in Theorem~1.2.1 in~\cite{sub}), the natural elliptic a priori estimate 
$$\exists C>0, \forall u \in D(q_1) \cap ... \cap D(q_N), \ \big\|\big(\langle(x,\xi)\rangle^{2}\big)^w u\big\|_{L^2} \leq C(\|q^w(x,D_x) u\|_{L^2}+\|u\|_{L^2}),$$
which easily implies (\ref{kristen5}).

We can therefore assume in the following that $k_0 \geq 1$ and find from Proposition~\ref{prop1} some real-valued weight functions 
\begin{equation}\label{giu8}\inc
g_j \in S\big(1,\langle X \rangle^{-\frac{2}{2k_0+1}}dX^2\big), \ 1 \leq j \leq N, \num
\end{equation}
such that 
\begin{equation}\label{giu9}\inc
\exists c,c_1,..., c_N>0, \forall X \in \rr^{2n}, \ 1+\sum_{j=1}^N\big(\textrm{Re }q_j(X)+c_jH_{\textrm{Im}q_j\ } g_j(X)\big) \geq c \langle X \rangle^{\frac{2}{2k_0+1}}. \num
\end{equation}
For $0<\eps \leq 1$, we consider the multipliers defined in the Wick quantization by symbols $1-\eps c_j g_j$. We recall that the definition of the Wick quantization and some elements of Wick calculus are recalled in Section~\ref{wick}. It follows from (\ref{giu8}), (\ref{lay0}), (\ref{lay1}), (\ref{lay2}) and the Cauchy-Schwarz inequality that 
\inc
\begin{multline*}\label{giu18}
\sum_{j=1}^N\textrm{Re}\big(q_j^{\textrm{Wick}} u, (1-\eps c_j g_j)^{\textrm{Wick}}u\big)= \sum_{j=1}^N\big(\textrm{Re}\big((1-\eps c_j g_j)^{\textrm{Wick}} q_j^{\textrm{Wick}}\big) u,u\big) \num \\
\leq  \sum_{j=1}^N \|1-\eps c_j g_j\|_{L^{\infty}}\|q_j^{\textrm{Wick}}u\|_{L^2}\|u\|_{L^2} 
\lesssim  \sum_{j=1}^N\|q_j^{\textrm{Wick}}u\|_{L^2}^2 +\|u\|_{L^2}^2 
\lesssim  \sum_{j=1}^N\|\tilde{q}_j^{w}u\|_{L^2}^2 +\|u\|_{L^2}^2,
\end{multline*}
where  
\begin{equation}\label{lay5}\inc
\tilde{q}_j(x,\xi)=q_j\Big(x,\frac{\xi}{2\pi}\Big),\num
\end{equation}
because the operators $(1-\eps c_j g_j)^{\textrm{Wick}}$ whose Wick symbol are real-valued, are formally selfadjoint. Indeed, symbols $r(q_j)$ defined in (\ref{lay2}) are here just some constants since $q_j$ are quadratic forms. The factor $2\pi$ in (\ref{lay5}) comes from the difference of normalizations chosen between (\ref{3}) and (\ref{lay3}) (See remark in Section~\ref{wick}). Since from (\ref{lay4}),
$$(1-\eps c_j g_j)^{\textrm{Wick}} q_j^{\textrm{Wick}} =\Big{[}(1-\eps c_j g_j)q_j+\frac{\eps}{4 \pi} c_j \nabla g_j. \nabla q_j-\frac{\eps}{4i \pi}c_j\{g_j,q_j\} \Big{]}^{\textrm{Wick}}+S_j,$$
with $\|S_j\|_{\mathcal{L}(L^2(\rr^n))} \lesssim 1$, we obtain from the fact that real Hamiltonians get quantized in the Wick quantization by formally selfadjoint operators that 
\begin{multline*}
\sum_{j=1}^N\textrm{Re}\big((1-\eps c_j g_j)^{\textrm{Wick}} q_j^{\textrm{Wick}}\big)= \sum_{j=1}^N\textrm{Re }S_j \\ +
\sum_{j=1}^N\Big{[}(1-\eps c_j g_j)\textrm{Re }q_j+\frac{\eps}{4 \pi}c_j \nabla g_j. \nabla \textrm{Re }q_j+\frac{\eps}{4 \pi}c_j H_{\textrm{Im}q_j}\ g_j \Big{]}^{\textrm{Wick}},
\end{multline*}
because $g_j$ are real-valued symbols. Since $\textrm{Re }q_j \geq 0$ and $g_j \in L^{\infty}(\rr^n)$, we can choose the positive parameter $\eps$ sufficiently small such that 
$$\forall \ 1 \leq j \leq N, \forall X \in \rr^{2n}, \ 1-\eps c_j g_j(X) \geq \frac{1}{2},$$
in order to deduce from (\ref{giu9}), (\ref{giu18}) and (\ref{lay0.5}) that 
\begin{equation}\label{lay11}\inc
\big((\langle X \rangle^{\frac{2}{2k_0+1}})^{\textrm{Wick}}u,u\big) \lesssim \|u\|_{L^2}^2+ \sum_{j=1}^N\|\tilde{q}_j^{w}u\|_{L^2}^2 + 
\sum_{j=1}^N\big|\big((\nabla g_j. \nabla \textrm{Re }q_j)^{\textrm{Wick}}u,u\big)\big|, \num
\end{equation}
because from (\ref{lay0.1}) and (\ref{lay0.2}), $1^{\textrm{Wick}}=\textrm{Id}.$

One can then complete the proof of Theorem~\ref{theorem1} by following exactly the same reasoning as the one used in \cite{sub}.  We recall this reasoning here for the sake of completeness of this work.

By denoting $\tilde{X}=\big(x,\xi/(2\pi)\big)$ and $\textrm{Op}^w\big(S(1,dX^2)\big)$ the operators obtained by the Weyl quantization of symbols in the class $S(1,dX^2)$, it follows from (\ref{lay1}), (\ref{lay2}) and usual results of symbolic calculus that 
\begin{equation}\label{lay43}\inc
\big(\langle X \rangle^{\frac{2}{2k_0+1}}\big)^{\textrm{Wick}} - \big(\langle \tilde{X} \rangle^{\frac{2}{2k_0+1}}\big)^{w} \in \textrm{Op}^w\big(S(1,dX^2)\big) \num
\end{equation}
and
\begin{equation}\label{lay44}\inc
\big(\langle \tilde{X} \rangle^{\frac{1}{2k_0+1}}\big)^{w}\big(\langle \tilde{X} \rangle^{\frac{1}{2k_0+1}}\big)^{w}- \big(\langle \tilde{X} \rangle^{\frac{2}{2k_0+1}}\big)^{w} \in \textrm{Op}^w\big(S(1,dX^2)\big),\num
\end{equation}
since $k_0 \geq 0$. By using that 
$$\big(\big(\langle \tilde{X} \rangle^{\frac{1}{2k_0+1}}\big)^{w}\big(\langle \tilde{X} \rangle^{\frac{1}{2k_0+1}}\big)^{w}u,u\big)=\big\|\big(\langle \tilde{X} \rangle^{\frac{1}{2k_0+1}}\big)^{w}u\big\|_{L^2}^2,$$
we therefore deduce from (\ref{lay11}) and the Calder\'on-Vaillancourt theorem that 
\begin{equation}\label{lay12}\inc
\big\|\big(\langle \tilde{X} \rangle^{\frac{1}{2k_0+1}}\big)^{w}u\big\|_{L^2}^2 \lesssim \|u\|_{L^2}^2+ \sum_{j=1}^N\|\tilde{q}_j^{w}u\|_{L^2}^2 + 
\sum_{j=1}^N\big|\big((\nabla g_j. \nabla \textrm{Re }q_j)^{\textrm{Wick}}u,u\big)\big|. \num
\end{equation}
Then, we get from (\ref{giu8}) and (\ref{lay0.5}) that 
\begin{equation}\label{lay13}\inc
\big|\big((\nabla g_j. \nabla \textrm{Re }q_j)^{\textrm{Wick}}u,u\big)\big| \lesssim \big(|\nabla \textrm{Re }q_j|^{\textrm{Wick}}u,u\big). \num
\end{equation}
Recalling now the well-known inequality 
\begin{equation}\label{giu00.1}\inc
|f'(x)|^2 \leq 2 f(x) \|f''\|_{L^{\infty}(\rr)}, \num
\end{equation}
fulfilled by any non-negative smooth function with bounded second derivative, we deduce from another use of (\ref{lay0.5}) that 
\begin{equation}\label{lay20}\inc
\big(|\nabla \textrm{Re }q_j|^{\textrm{Wick}}u,u\big) \lesssim \big(((\textrm{Re }q_j)^{\frac{1}{2}})^{\textrm{Wick}}u,u\big) \lesssim \big((1+ \textrm{Re }q_j)^{\textrm{Wick}}u,u\big), \num
\end{equation}
since $\textrm{Re }q_j$ is a non-negative quadratic form and that 
$$2(\textrm{Re }q_j)^{\frac{1}{2}} \leq 1+\textrm{Re }q_j.$$ 
By using the same arguments as in (\ref{giu18}), we obtain that 
\begin{multline*}
\big((1+ \textrm{Re }q_j)^{\textrm{Wick}}u,u\big)=\big((\textrm{Re }q_j)^{\textrm{Wick}}u,u\big)+\|u\|_{L^2}^2=\textrm{Re}(q_j^{\textrm{Wick}}u,u)+\|u\|_{L^2}^2\\
\leq \|q_j^{\textrm{Wick}}u\|_{L^2}\|u\|_{L^2}+\|u\|_{L^2}^2
\lesssim  \|q_j^{\textrm{Wick}}u\|_{L^2}^2 +\|u\|_{L^2}^2 \lesssim  \|\tilde{q}_j^{w}u\|_{L^2}^2 +\|u\|_{L^2}^2.
\end{multline*}
It therefore follows from (\ref{lay12}), (\ref{lay13}) and (\ref{lay20}) that 
\begin{equation}\label{lay21}\inc
\big\|\big(\langle \tilde{X} \rangle^{\frac{1}{2k_0+1}}\big)^{w}u\big\|_{L^2}^2 \lesssim  \|u\|_{L^2}^2+ \sum_{j=1}^N\|\tilde{q}_j^{w}u\|_{L^2}^2 . \num
\end{equation}
In order to improve the estimate (\ref{lay21}), we carefully resume our previous analysis and notice that our previous reasoning has in fact established that
\begin{align*}
& \ \big\|\big(\langle \tilde{X} \rangle^{\frac{1}{2k_0+1}}\big)^{w}u\big\|_{L^2}^2 \\
\lesssim & \  \|u\|_{L^2}^2+\sum_{j=1}^N\big|\textrm{Re}\big(q_j^{\textrm{Wick}} u, (1-\eps c_j g_j)^{\textrm{Wick}}u\big)\big|+ 
\sum_{j=1}^N\big|\big((\nabla g_j. \nabla \textrm{Re }q_j)^{\textrm{Wick}}u,u\big)\big| \\
\lesssim & \  \|u\|_{L^2}^2+ \sum_{j=1}^N\big|\textrm{Re}\big(q_j^{\textrm{Wick}} u, (1-\eps c_j g_j)^{\textrm{Wick}}u\big)\big|+ 
\sum_{j=1}^N|\textrm{Re}(q_j^{\textrm{Wick}} u, u)|\\
\lesssim & \  \|u\|_{L^2}^2+\sum_{j=1}^N\big|\textrm{Re}\big(\tilde{q}_j^{w} u, (1-\eps c_j g_j)^{\textrm{Wick}}u\big)\big|+ 
\sum_{j=1}^N|\textrm{Re}(\tilde{q}_j^{w} u, u)|, 
\end{align*}
because $(1-\eps c_j g_j)^{\textrm{Wick}}$ is a bounded operator on $L^2(\rr^n)$,
\begin{equation}\label{linn1}\inc
\|(1-\eps c_j g_j)^{\textrm{Wick}}\|_{\mathcal{L}(L^2)} \leq \|1-\eps c_j g_j\|_{L^{\infty}(\rr^{2n})}.\num
\end{equation}
By applying this estimate to $\big(\langle \tilde{X} \rangle^{\frac{1}{2k_0+1}}\big)^{w}u$, we deduce from (\ref{lay44}) and the Calder\'on-Vaillancourt theorem that 
\inc\begin{multline*}\label{lay45}
\big\|\big(\langle \tilde{X} \rangle^{\frac{2}{2k_0+1}}\big)^{w}u\big\|_{L^2}^2 \lesssim   
\sum_{j=1}^N\Big|\textrm{Re}\Big(\tilde{q}_j^{w} \big(\langle \tilde{X} \rangle^{\frac{1}{2k_0+1}}\big)^{w}u, \big(\langle \tilde{X} \rangle^{\frac{1}{2k_0+1}}\big)^{w}u\Big)\Big|\\
+\sum_{j=1}^N\Big|\textrm{Re}\Big(\tilde{q}_j^{w} \big(\langle \tilde{X} \rangle^{\frac{1}{2k_0+1}}\big)^{w} u, 
(1-\eps c_j g_j)^{\textrm{Wick}}\big(\langle \tilde{X} \rangle^{\frac{1}{2k_0+1}}\big)^{w}u\Big)\Big|+\big\|\big(\langle \tilde{X} \rangle^{\frac{1}{2k_0+1}}\big)^{w}u\big\|_{L^2}^2+\|u\|_{L^2}^2. \num
\end{multline*}
Then, by noticing that the commutator 
\begin{equation}\label{lay50}\inc
\big[\tilde{q}_j^w,\big(\langle \tilde{X} \rangle^{\frac{1}{2k_0+1}}\big)^{w}\big] \in \textrm{Op}^w\big(S\big(\langle X \rangle^{\frac{1}{2k_0+1}},\langle X \rangle^{-2}dX^2\big)\big), \num
\end{equation}
because $\tilde{q}_j$ is a quadratic form, and that 
\begin{equation}\label{lay51}\inc
\big(\langle \tilde{X} \rangle^{-\frac{1}{2k_0+1}}\big)^{w}\big(\langle \tilde{X} \rangle^{\frac{1}{2k_0+1}}\big)^{w}- \textrm{Id} \in \textrm{Op}^w\big(S(\langle X \rangle^{-2},\langle X \rangle^{-2}dX^2)\big),\num
\end{equation}
we deduce from standard results of symbolic calculus and the Calder\'on-Vaillancourt theorem that 
\inc\begin{align*}\label{lay52}
\big\|\big[\tilde{q}_j^w,\big(\langle \tilde{X} \rangle^{\frac{1}{2k_0+1}}\big)^{w}\big]u\big\|_{L^2} \lesssim & \ \big\|\big[\tilde{q}_j^w,\big(\langle \tilde{X} \rangle^{\frac{1}{2k_0+1}}\big)^{w}\big]\big(\langle \tilde{X} \rangle^{-\frac{1}{2k_0+1}}\big)^{w}\big(\langle \tilde{X} \rangle^{\frac{1}{2k_0+1}}\big)^{w}u\big\|_{L^2}+ \|u\|_{L^2}\\
\lesssim & \ \big\|\big(\langle \tilde{X} \rangle^{\frac{1}{2k_0+1}}\big)^{w}u\big\|_{L^2}+\|u\|_{L^2}. \num
\end{align*}
By introducing this commutator, we get from
the Cauchy-Schwarz inequality and (\ref{lay52}) that 
\begin{multline*}
\Big|\textrm{Re}\Big(\tilde{q}_j^{w} \big(\langle \tilde{X} \rangle^{\frac{1}{2k_0+1}}\big)^{w}u, \big(\langle \tilde{X} \rangle^{\frac{1}{2k_0+1}}\big)^{w}u\Big)\Big| \lesssim 
\Big|\textrm{Re}\Big(\tilde{q}_j^{w} u, \big(\langle \tilde{X} \rangle^{\frac{1}{2k_0+1}}\big)^{w}\big(\langle \tilde{X} \rangle^{\frac{1}{2k_0+1}}\big)^{w}u\Big)\Big|
\\
+\big\|\big(\langle \tilde{X} \rangle^{\frac{1}{2k_0+1}}\big)^{w}u\big\|_{L^2}^2+\|u\|_{L^2}^2.
\end{multline*}
Another use of the Cauchy-Schwarz inequality and the Calder\'on-Vaillancourt theorem with (\ref{lay44}) gives that 
$$\Big|\textrm{Re}\Big(\tilde{q}_j^{w} u, \big(\langle \tilde{X} \rangle^{\frac{1}{2k_0+1}}\big)^{w}\big(\langle \tilde{X} \rangle^{\frac{1}{2k_0+1}}\big)^{w}u\Big)\Big| \lesssim
\|\tilde{q}_j^{w} u\|_{L^2}\big\|\big(\langle \tilde{X} \rangle^{\frac{2}{2k_0+1}}\big)^{w}u\big\|_{L^2}+\|\tilde{q}_j^{w} u\|_{L^2}\|u\|_{L^2}.$$
We then deduce from (\ref{lay21}) and the previous estimate that
\begin{align*}
& \ \sum_{j=1}^N\Big|\textrm{Re}\Big(\tilde{q}_j^{w} \big(\langle \tilde{X} \rangle^{\frac{1}{2k_0+1}}\big)^{w}u, \big(\langle \tilde{X} \rangle^{\frac{1}{2k_0+1}}\big)^{w}u\Big)\Big| \\ 
\lesssim & \
\big\|\big(\langle \tilde{X} \rangle^{\frac{2}{2k_0+1}}\big)^{w}u\big\|_{L^2} \sum_{j=1}^N \|\tilde{q}_j^{w} u\|_{L^2}
+\sum_{j=1}^N\|\tilde{q}_j^{w} u\|_{L^2}^2+\|u\|_{L^2}^2.
\end{align*}
By using again the Cauchy-Schwarz inequality, (\ref{lay21}), (\ref{linn1}), (\ref{lay45}) and (\ref{lay52}), this estimate implies that 
\inc\begin{align*}\label{lay46}
& \ \big\|\big(\langle \tilde{X} \rangle^{\frac{2}{2k_0+1}}\big)^{w}u\big\|_{L^2}^2 \lesssim  \sum_{j=1}^N\Big|\textrm{Re}\Big(\big[\tilde{q}_j^{w}, \big(\langle \tilde{X} \rangle^{\frac{1}{2k_0+1}}\big)^{w}\big] u, 
(1-\eps c_j g_j)^{\textrm{Wick}}\big(\langle \tilde{X} \rangle^{\frac{1}{2k_0+1}}\big)^{w}u\Big)\Big| \num \\ 
+& \ \sum_{j=1}^N\Big|\textrm{Re}\Big(\tilde{q}_j^{w}  u, 
\big(\langle \tilde{X} \rangle^{\frac{1}{2k_0+1}}\big)^{w}(1-\eps c_j g_j)^{\textrm{Wick}}\big(\langle \tilde{X} \rangle^{\frac{1}{2k_0+1}}\big)^{w}u\Big)\Big|
+\sum_{j=1}^N\|\tilde{q}_j^{w} u\|_{L^2}^2+\|u\|_{L^2}^2 \\
\lesssim & \ \sum_{j=1}^N\Big|\textrm{Re}\Big(\tilde{q}_j^{w}  u, 
\big(\langle \tilde{X} \rangle^{\frac{1}{2k_0+1}}\big)^{w}(1-\eps c_j g_j)^{\textrm{Wick}}\big(\langle \tilde{X} \rangle^{\frac{1}{2k_0+1}}\big)^{w}u\Big)\Big|
+\sum_{j=1}^N\|\tilde{q}_j^{w} u\|_{L^2}^2+\|u\|_{L^2}^2\\
\lesssim & \ \sum_{j=1}^N\|\tilde{q}_j^{w}  u\|_{L^2}\big\| 
\big(\langle \tilde{X} \rangle^{\frac{1}{2k_0+1}}\big)^{w}(1-\eps c_j g_j)^{\textrm{Wick}}\big(\langle \tilde{X} \rangle^{\frac{1}{2k_0+1}}\big)^{w}u\big\|_{L^2}
+\sum_{j=1}^N\|\tilde{q}_j^{w} u\|_{L^2}^2+\|u\|_{L^2}^2,
\end{align*}
because we get from (\ref{linn1}) and (\ref{lay52}) that 
\begin{multline*}
\Big|\textrm{Re}\Big(\big[\tilde{q}_j^{w}, \big(\langle \tilde{X} \rangle^{\frac{1}{2k_0+1}}\big)^{w}\big] u, 
(1-\eps c_j g_j)^{\textrm{Wick}}\big(\langle \tilde{X} \rangle^{\frac{1}{2k_0+1}}\big)^{w}u\Big)\Big| \lesssim 
\big\|\big(\langle \tilde{X} \rangle^{\frac{1}{2k_0+1}}\big)^{w}u\big\|_{L^2}^2\\ +\big\|\big(\langle \tilde{X} \rangle^{\frac{1}{2k_0+1}}\big)^{w}u\big\|_{L^2}\|u\|_{L^2}.
\end{multline*}
Notice now that (\ref{giu8}), (\ref{lay1bis}) and (\ref{lay2bis}) imply that 
$$\big[\big(\langle \tilde{X} \rangle^{\frac{1}{2k_0+1}}\big)^{w},(1-\eps c_j g_j)^{\textrm{Wick}}\big] \in \textrm{Op}^w\big(S(1,dX^2)\big),$$
since $(1-\eps c_j g_j)^{\textrm{Wick}}=\tilde{g}_j^w$, with $\tilde{g}_j \in S(1,dX^2)$ and $k_0 \geq 0$. 
By introducing this new commutator, we deduce from the Calder\'on-Vaillancourt theorem, (\ref{lay44}), (\ref{lay21}) and (\ref{linn1}) that
\begin{align*}
& \ \big\| \big(\langle \tilde{X} \rangle^{\frac{1}{2k_0+1}}\big)^{w}(1-\eps c_j g_j)^{\textrm{Wick}}\big(\langle \tilde{X} \rangle^{\frac{1}{2k_0+1}}\big)^{w}u\big\|_{L^2} \\
\lesssim & \ \big\|\big(\langle \tilde{X} \rangle^{\frac{1}{2k_0+1}}\big)^{w}u\big\|_{L^2}
+\big\| (1-\eps c_j g_j)^{\textrm{Wick}}\big(\langle \tilde{X} \rangle^{\frac{1}{2k_0+1}}\big)^{w}\big(\langle \tilde{X} \rangle^{\frac{1}{2k_0+1}}\big)^{w}u\big\|_{L^2} \\
\lesssim & \ \big\|\big(\langle \tilde{X} \rangle^{\frac{1}{2k_0+1}}\big)^{w}u\big\|_{L^2}
+\big\|\big(\langle \tilde{X} \rangle^{\frac{1}{2k_0+1}}\big)^{w}\big(\langle \tilde{X} \rangle^{\frac{1}{2k_0+1}}\big)^{w}u\big\|_{L^2} \\\lesssim & \ \big\|\big(\langle \tilde{X} \rangle^{\frac{2}{2k_0+1}}\big)^{w}u\big\|_{L^2}+\big\|\big(\langle \tilde{X} \rangle^{\frac{1}{2k_0+1}}\big)^{w}u\big\|_{L^2}+\|u\|_{L^2}\\
\lesssim & \ \big\|\big(\langle \tilde{X} \rangle^{\frac{2}{2k_0+1}}\big)^{w}u\big\|_{L^2}+\sum_{j=1}^{N}\|\tilde{q}_j^{w}u\|_{L^2}+\|u\|_{L^2}.
\end{align*}
Recalling (\ref{lay46}), we can then use this last estimate to obtain that 
\begin{equation}\label{lay60}\inc
\big\|\big(\langle \tilde{X} \rangle^{\frac{2}{2k_0+1}}\big)^{w}u\big\|_{L^2}^2 \lesssim \sum_{j=1}^N\|\tilde{q}_j^{w}u\|_{L^2}^2+\|u\|_{L^2}^2. \num
\end{equation}
By finally noticing from the homogeneity of degree 2 of $\tilde{q}_j$ that we have 
$$(\tilde{q}_j \circ T)(x,\xi)=\frac{1}{2\pi}q_j(x,\xi),$$
if $T$ stands for the real linear symplectic transformation 
$$T(x,\xi)=\big((2\pi)^{-\frac{1}{2}}x,(2\pi)^{\frac{1}{2}}\xi\big),$$
we deduce from the symplectic invariance of the Weyl quantization (Theorem~18.5.9 in~\cite{hormander}) that  
$$\big\|\big(\langle X \rangle^{\frac{2}{2k_0+1}}\big)^{w}u\big\|_{L^2}^2 \lesssim  \sum_{j=1}^N\|q_j^{w}u\|_{L^2}^2 +\|u\|_{L^2}^2,$$
which proves Theorem~\ref{theorem1}.

\section{Proof of Proposition~\ref{prop1}}\label{proofprop1}
\init

We prove Proposition~\ref{prop1} by induction on the positive integer $m \geq 1$ appearing in (\ref{giu2}). Let $m \geq 1$,  we shall assume that Proposition~\ref{prop1} is fulfilled for any 
open set $\Omega_0$ of $\rr^{2n}$, when the positive integer in (\ref{giu2}) is strictly smaller than $m$.

In the following, we denote by $\psi$, $\chi$ and $w$ some $C^{\infty}(\rr,[0,1])$ functions respectively satisfying  
\begin{equation}\label{giu13}\inc
\psi=1 \textrm{ on } [-1,1], \ \textrm{supp } \psi \subset [-2,2], \num
\end{equation} 
\begin{equation}\label{giu14}\inc
\chi=1 \textrm{ on } \{x \in \rr : 1 \leq |x| \leq 2 \}, \ \textrm{supp } \chi \subset \{x \in \rr :  1/2 \leq |x| \leq 3\}, \num
\end{equation}
and
\begin{equation}\label{giu15}\inc
w=1 \textrm{ on } \{x \in \rr : |x| \geq 2 \}, \ \textrm{supp } w \subset \{x \in \rr :  |x| \geq 1\}. \num
\end{equation}
More generically, we shall denote by $\psi_j$, $\chi_j$ and $w_j$, $j \in \nn$, some other $C^{\infty}(\rr,[0,1])$ functions satisfying similar properties as respectively $\psi$, $\chi$ and $w$ with possibly different choices for the positive numerical values which define their support localizations.

Let $\Omega_0$ be an open set of $\rr^{2n}$ such that (\ref{giu2}) is fulfilled. Considering the quadratic forms 
\begin{equation}\label{giu11}\inc
\tilde{r}_{1,p}(X)=\sum_{j=1}^N\textrm{Re }q_j(X;\textrm{Im }F_{p}X),  \num
\end{equation}
\begin{equation}\label{giu11serena}\inc
\tilde{r}_{k,p}(X)=\hspace{-0.8cm}\sum_{\substack{j=1,...,N \\ (l_1,...,l_{k-1}) \in \{1,...,N\}^{k-1}}}\hspace{-0.8cm} \textrm{Re }q_j(\textrm{Im }F_{l_1}...\textrm{Im }F_{l_{k-1}}X;\textrm{Im }F_{l_1}...\textrm{Im }F_{l_{k-1}}\textrm{Im }F_{p}X),  \num
\end{equation}
for any $1 \leq p \leq N$, $2 \leq k \leq m$; 
\begin{equation}\label{bellagiu11}\inc
r_0(X)=\sum_{j=1}^N\textrm{Re }q_j(X), \qquad r_{k}(X)=\hspace{-0.5cm}\sum_{\substack{j=1,...,N \\ (l_1,...,l_{k}) \in \{1,...,N\}^k}}\hspace{-0.5cm} \textrm{Re }q_j(\textrm{Im }F_{l_1}...\textrm{Im }F_{l_{k}}X),  \num
\end{equation}
for any $1 \leq k \leq m$; and defining
\begin{equation}\label{giu23}\inc
\tilde{g}_{m,p}(X)=\psi\big(r_{m-1}(X)\langle X \rangle^{-\frac{2(2m-1)}{2m+1}}\big)\langle X \rangle^{-\frac{4m}{2m+1}}\tilde{r}_{m,p}(X), \num
\end{equation}
where $\psi$ is the function defined in (\ref{giu13}) and $1 \leq p \leq N$, we get from Lemma~\ref{lem2} that 
\begin{align*}\label{giu24}\inc
& \ H_{\textrm{Im}q_p}\ \tilde{g}_{m,p}(X)=2\psi\big(r_{m-1}(X)\langle X \rangle^{-\frac{2(2m-1)}{2m+1}}\big)\hspace{-1.3cm} \sum_{\substack{j=1,...,N \\ (l_1,...,l_{m-1}) \in \{1,...,N\}^{m-1}}}\hspace{-1.3cm}\frac{\textrm{Re }q_j(\textrm{Im }F_{l_1}...\textrm{Im }F_{l_{m-1}}\textrm{Im }F_{p}X)}{\langle X \rangle^{\frac{4m}{2m+1}}} \num \\
+ & \ 2\psi\big(r_{m-1}(X)\langle X \rangle^{-\frac{2(2m-1)}{2m+1}}\big)\hspace{-1.4cm} \sum_{\substack{j=1,...,N \\ (l_1,...,l_{m-1}) \in \{1,...,N\}^{m-1}}} \hspace{-1.4cm} \frac{ \textrm{Re }q_j(\textrm{Im }F_{l_1}...\textrm{Im }F_{l_{m-1}}X;\textrm{Im }F_{l_1}...\textrm{Im }F_{l_{m-1}}(\textrm{Im }F_{p})^2X)}{\langle X \rangle^{\frac{4m}{2m+1}}}\\
+ & \ H_{\textrm{Im}q_p}\Big(\psi\big(r_{m-1}(X) \langle X \rangle^{-\frac{2(2m-1)}{2m+1}}\big)\Big)\frac{\tilde{r}_{m,p}(X)}{\langle X \rangle^{\frac{4m}{2m+1}}}\\
+ & \ \psi\big(r_{m-1}(X) \langle X \rangle^{-\frac{2(2m-1)}{2m+1}}\big)H_{\textrm{Im}q_p}\big(\langle X \rangle^{-\frac{4m}{2m+1}}\big)\tilde{r}_{m,p}(X).
\end{align*}
We first check that 
\begin{equation}\label{giu25}\inc
\tilde{g}_{m,p} \in S\big(1,\langle X \rangle^{-\frac{2(2m-1)}{2m+1}}dX^2\big).\num
\end{equation}
In order to verify this, we notice from Lemma~\ref{lem2.21} that the quadratic forms 
\begin{equation}\label{lay80}\inc
\textrm{Re }q_j(\textrm{Im }F_{l_1}...\textrm{Im }F_{l_{m-1}}X;\textrm{Im }F_{l_1}...\textrm{Im }F_{l_{m-1}}\textrm{Im }F_{p}X)\num
\end{equation} 
and
\begin{equation}\label{lay80fav}\inc
 \textrm{Re }q_j(\textrm{Im }F_{l_1}...\textrm{Im }F_{l_{m-1}}X;\textrm{Im }F_{l_1}...\textrm{Im }F_{l_{m-1}}(\textrm{Im }F_{p})^2X),\num
\end{equation} 
belong to the symbol class 
\begin{equation}\label{lay81}\inc
S_{\Omega}\big(\langle X \rangle^{\frac{4m}{2m+1}}, \langle X \rangle^{-\frac{2(2m-1)}{2m+1}}dX^2\big),\num
\end{equation}
for any open set $\Omega$ in $\rr^{2n}$ where $r_{m-1}(X) \lesssim \langle X \rangle^{\frac{2(2m-1)}{2m+1}}$. To check this, we just use in addition to Lemma~\ref{lem2.21} the obvious estimates
$$\textrm{Re }q_j(\textrm{Im }F_{l_1}...\textrm{Im }F_{l_{m-1}}\textrm{Im }F_{p}X)^{\frac{1}{2}} \lesssim \langle X \rangle$$
and  
$$\textrm{Re }q_j(\textrm{Im }F_{l_1}...\textrm{Im }F_{l_{m-1}}(\textrm{Im }F_{p})^2X)^{\frac{1}{2}} \lesssim \langle X \rangle.$$
Moreover, since 
\begin{equation}\label{oc1}\inc
\langle X \rangle^{-\frac{4m}{2m+1}} \in S\big(\langle X \rangle^{-\frac{4m}{2m+1}},\langle X \rangle^{-2}dX^2\big),\num
\end{equation}
we obtain (\ref{giu25}) from (\ref{giu13}), (\ref{giu11serena}), (\ref{bellagiu11}), (\ref{giu23}), (\ref{lay80}), (\ref{lay81}) and Lemma~\ref{lem3}.

Denoting respectively $A_{1,p}$, $A_{2,p}$, $A_{3,p}$ and $A_{4,p}$ the four terms appearing in the right hand side of (\ref{giu24}), we first notice from (\ref{giu13}), (\ref{lay80}), (\ref{lay81}), (\ref{oc1}) and Lemma~\ref{lem3} that 
\begin{equation}\label{giu27}\inc
A_{2,p} \in  S\big(1,\langle X \rangle^{-\frac{2(2m-1)}{2m+1}}dX^2\big).\num
\end{equation}
Next, by using that 
$$\textrm{Im }q_p \in S\big(\langle X \rangle^2,\langle X \rangle^{-2}dX^2\big),$$
since $\textrm{Im }q_p$ is a quadratic form,
we get from (\ref{giu13}), (\ref{giu11serena}), (\ref{bellagiu11}), (\ref{lay80}), (\ref{lay81}), (\ref{oc1}) and Lemma~\ref{lem3} that
\begin{equation}\label{lay82}\inc
A_{3,p} \in S\big(\langle X \rangle^{\frac{2}{2m+1}},\langle X \rangle^{-\frac{2(2m-1)}{2m+1}}dX^2\big),\num
\end{equation}
since
$$H_{\textrm{Im}q_p}\Big(\psi\big(r_{m-1}(X)\langle X \rangle^{-\frac{2(2m-1)}{2m+1}}\big)\Big)\in S\big(\langle X \rangle^{\frac{2}{2m+1}},\langle X \rangle^{-\frac{2(2m-1)}{2m+1}}dX^2\big).$$
By using now that 
$$H_{\textrm{Im}q_p}\big(\langle X \rangle^{-\frac{4m}{2m+1}}\big) \in S\big(\langle X \rangle^{-\frac{4m}{2m+1}},\langle X \rangle^{-2}dX^2\big),$$
we finally obtain from another use of (\ref{giu13}), (\ref{giu11serena}), (\ref{bellagiu11}), (\ref{lay80}), (\ref{lay81}) and Lemma~\ref{lem3} that 
\begin{equation}\label{lay83}\inc
A_{4,p} \in S\big(1,\langle X \rangle^{-\frac{2(2m-1)}{2m+1}}dX^2\big).\num
\end{equation}
Since the term $A_{3,p}$ is supported in 
$$\textrm{supp } \psi'\big(r_{m-1}(X)\langle X \rangle^{-\frac{2(2m-1)}{2m+1}}\big),$$ 
we deduce from (\ref{giu24}), (\ref{giu27}), (\ref{lay82}) and (\ref{lay83}) that there exists $\chi_0$ a $C^{\infty}(\rr,[0,1])$ function satisfying similar properties as in (\ref{giu14}), with possibly different positive numerical values for its support localization, such that, $\exists c_{1},c_{2}>0$, $\forall X \in \rr^{2n}$,
\begin{align*}\label{giu28}\inc
 & \ c_{1}+
c_{2}\chi_0\big(r_{m-1}(X)\langle X \rangle^{-\frac{2(2m-1)}{2m+1}}\big)\langle X \rangle^{\frac{2}{2m+1}}+\sum_{p=1}^NH_{\textrm{Im}q_p}\ \tilde{g}_{m,p}(X) \num \\
 & \ \geq 2 \psi\big(r_{m-1}(X)\langle X \rangle^{-\frac{2(2m-1)}{2m+1}}\big)  \frac{r_m(X)}{\langle X \rangle^{\frac{4m}{2m+1}}}. 
\end{align*}
Recalling (\ref{giu2}), one can find some positive constants $c_3, c_4>0$ such that
\begin{equation}\label{giu29}\inc  
 \sum_{k=0}^{m-1}r_k(X)  \geq c_3 |X|^2, \num
\end{equation}on the open set 
\begin{equation}\label{oc2}\inc
\Omega_1=\left\{X \in \rr^{2n} : r_{m}(X) < c_4|X|^2\right\} \cap \Omega_0. \num
\end{equation}
When $m \geq 2$, one can find according to our induction hypothesis some real-valued functions 
\begin{equation}\label{giu29.1}\inc
\tilde{\mathfrak{g}}_{m,p} \in S_{\Omega_1}\big(1,\langle X \rangle^{-\frac{2}{2m-1}}dX^2\big), \ 1 \leq p \leq N, \num
\end{equation}
such that 
\begin{equation}\label{giu30}\inc
\exists c_{5,p}>0, \forall X \in \Omega_1, \ 1+\sum_{p=1}^N\big(\textrm{Re }q_p(X)+c_{5,p}H_{\textrm{Im}q_p\ } \tilde{\mathfrak{g}}_{m,p}(X)\big)  \gtrsim  \langle X \rangle^{\frac{2}{2m-1}}. \num
\end{equation}
For convenience, we set in the following $\tilde{\mathfrak{g}}_{1,p}=0$ when $m=1$.  
By choosing suitably $\psi_0$ and $w_0$ some $C^{\infty}(\rr,[0,1])$ functions satisfying similar properties as the functions respectively defined in (\ref{giu13}) and (\ref{giu15}), with possibly different positive numerical values for their support localizations, such that 
\inc\begin{equation}\label{bri2}
\textrm{supp } \psi_0\big(r_m(X)|X|^{-2}\big)w_0(X)  \subset \big\{X \in \rr^{2n} : r_m(X) < c_4|X|^2\big\},\num
\end{equation}
and setting 
\begin{equation}\label{giu32}\inc
G_{m,p}(X)=\tilde{g}_{m,p}(X)+ \psi_0\big(r_m(X) |X|^{-2}\big)w_0(X)\tilde{\mathfrak{g}}_{m,p}(X), \ X \in \Omega_0, \num
\end{equation} 
we deduce from a straightforward adaptation of the Lemma~\ref{lem3} by recalling (\ref{giu13}) and (\ref{giu15}) that 
\begin{equation}\label{bri1}\inc
\psi_0\big(r_m(X)|X|^{-2}\big)w_0(X) \in S\big(1,\langle X \rangle^{-2}dX^2\big).\num
\end{equation}
According to (\ref{giu25}) and (\ref{giu29.1}), this implies that
\begin{equation}\label{giu33}\inc
G_{1,p} \in S_{\Omega_0}\big(1,\langle X \rangle^{-\frac{2}{3}}dX^2\big) \textrm{ and } G_{m,p} \in S_{\Omega_0}\big(1,\langle X \rangle^{-\frac{2}{2m-1}}dX^2\big), \num
\end{equation}
when $m \geq 2$.
Since from (\ref{bri1}),
$$H_{\textrm{Im}q_p}\big(\psi_0\big(r_m(X)|X|^{-2}\big)w_0(X) \big) \in S\big(1,\langle X \rangle^{-2}dX^2\big),$$
because $\textrm{Im }q_p$ is a quadratic form, 
we first notice from (\ref{oc2}), (\ref{giu29.1}) and (\ref{bri2}) that 
$$H_{\textrm{Im}q_p}\big(\psi_0\big(r_m(X)|X|^{-2}\big)w_0(X) \big)\tilde{\mathfrak{g}}_{m,p}(X)
\in S_{\Omega_0}\big(1,\langle X \rangle^{-\frac{2}{2m-1}}dX^2\big),$$ 
and then deduce from (\ref{giu28}), (\ref{oc2}), (\ref{giu30}), (\ref{bri2}) and (\ref{giu32}) that there exist some positive contants $c_{6,p},c_7>0$ such that for all $X \in \Omega_0$,
\begin{align*}
& \ \sum_{p=1}^N\big(\textrm{Re }q_p(X)+c_{6,p} H_{\textrm{Im}q_p}\ G_{m,p}(X)\big)+1+
c_7\chi_0\big(r_{m-1}(X)\langle X \rangle^{-\frac{2(2m-1)}{2m+1}}\big)\langle X \rangle^{\frac{2}{2m+1}}  \\
& \ \gtrsim \psi\big(r_{m-1}(X)\langle X \rangle^{-\frac{2(2m-1)}{2m+1}}\big)\frac{r_m(X)}{\langle X \rangle^{\frac{4m}{2m+1}}}
  +\psi_0\big(r_{m}(X)|X |^{-2}\big)w_0(X)\langle X \rangle^{\frac{2}{2m-1}}, 
\end{align*}
when $m \geq 2$.
Since 
$$\langle X \rangle^{\frac{2}{2m-1}} \gtrsim \langle X \rangle^{\frac{2}{2m+1}} \textrm{ and }
\frac{r_m(X)}{\langle X \rangle^{\frac{4m}{2m+1}}} \gtrsim |X|^{\frac{2}{2m+1}},$$
when $r_m(X) \gtrsim  |X|^2$, 
we deduce from the previous estimate by distinguishing the regions in $\Omega_0$ where
$$r_m(X) \lesssim  |X|^2 \textrm{ and } r_m(X) \gtrsim  |X|^2,$$ 
according to the support of the function
$$\psi_0\big(r_{m}(X)|X |^{-2}\big),$$
that one can find  a $C^{\infty}(\rr,[0,1])$ function $w_1$ with the same kind of support as the function defined in (\ref{giu15}) such that 
\inc\begin{multline*}\label{giu35}
\exists c_{8,p},c_9>0, \forall X \in \Omega_0, \  \sum_{p=1}^N\big(\textrm{Re }q_p(X) \num +c_{8,p}H_{\textrm{Im}q_p}\ G_{m,p}(X)\big)\\ 
+c_9 w_1\big(r_{m-1}(X)\langle X \rangle^{-\frac{2(2m-1)}{2m+1}}\big)\langle X \rangle^{\frac{2}{2m+1}} 
+1 \gtrsim \langle X \rangle^{\frac{2}{2m+1}}, 
\end{multline*}
when $m \geq 2$. When $m=1$, 
we notice from (\ref{giu2}) that
\begin{equation}\label{bri3}\inc
r_1(X) \gtrsim \langle X \rangle^{2},\num
\end{equation}
on any set where
\begin{equation}\label{bri4}\inc
|X|\geq c_{10} \textrm{ and } r_0(X)=\sum_{p=1}^N\textrm{Re }q_p(X) \leq  \langle X \rangle^{\frac{2}{3}},\num
\end{equation}
if the positive constant $c_{10}$ is chosen sufficiently large.  Moreover, since in this case $G_{1,p}=\tilde{g}_{1,p}$ and that $\textrm{Re }q_p \geq 0$, one can deduce from (\ref{giu13}), (\ref{giu15}), (\ref{giu28}), (\ref{bri3}) and (\ref{bri4}),
by distinguishing the regions in $\Omega_0$ where 
$$r_0(X) \lesssim \langle X \rangle^{\frac{2}{3}} \textrm{ and } r_0(X) \gtrsim \langle X \rangle^{\frac{2}{3}},$$
according to the support of the function
$$\psi\big(r_0(X)\langle X \rangle^{-\frac{2}{3}}\big),$$
that the estimate (\ref{giu35}) is also fulfilled in the case $m=1$. Continuing our study of the case where $m=1$, we notice from (\ref{giu15}) and $\textrm{Re }q_p \geq 0$, that one can estimate
$$w_1\big(r_0(X)\langle X \rangle^{-\frac{2}{3}}\big)\langle X \rangle^{\frac{2}{3}} \lesssim r_0(X)=\sum_{p=1}^N\textrm{Re }q_p(X),$$
for all $X \in \rr^{2n}$. It therefore follows that one can find $c_{11,p}>0$ such that for all $X \in \Omega_0$,
$$\sum_{p=1}^N\big(\textrm{Re }q_p(X) +c_{11,p}H_{\textrm{Im}q_p} G_{1,p}(X)\big) +1 \gtrsim \langle X \rangle^{\frac{2}{3}},$$
which proves Proposition~\ref{prop1} in the case where $m=1$, and our induction hypothesis in the basis case.

Assuming in the following that $m \geq 2$, we shall now work on the term 
$$w_1\big(r_{m-1}(X)\langle X \rangle^{-\frac{2(2m-1)}{2m+1}}\big)\langle X \rangle^{\frac{2}{2m+1}},$$
appearing in (\ref{giu35}).
By considering some constants $\Lambda_j \geq 1$, for $0 \leq j \leq m-2$, whose values will be successively chosen in the following, we shall prove that one can write that for all $X \in \rr^{2n}$,
\inc\begin{multline*}\label{giu38}
w_1\left(\frac{r_{m-1}(X)}{\langle X \rangle^{\frac{2(2m-1)}{2m+1}}}\right) \leq \tilde{W}_0(X)\Psi_0(X) \\ + \sum_{j=1}^{m-2}{\tilde{W}_0(X)
\Big(\prod_{l=1}^{j}W_l(X)\Big)\Psi_j(X)}
+\tilde{W}_{0}(X)\Big(\prod_{l=1}^{m-1}W_l(X)\Big), \num
\end{multline*}
with 
\begin{equation}\label{giu39}\inc
\Psi_j(X)=\psi\left(\frac{\Lambda_j r_{m-j-2}(X)}{r_{m-j-1}(X)^{\frac{2m-2j-3}{2m-2j-1}}}\right), \ 0 \leq j \leq m-2, \num
\end{equation} 
\begin{equation}\label{giu40}\inc
W_j(X)=w_{2}\left(\frac{\Lambda_{j-1} r_{m-j-1}(X)}{r_{m-j}(X)^{\frac{2m-2j-1}{2m-2j+1}}}\right), \ 1 \leq j \leq m-1,\num
\end{equation}
\begin{equation}\label{giu41}\inc
\tilde{W}_0(X)=w_{1}\left(\frac{r_{m-1}(X)}{\langle X \rangle^{\frac{2(2m-1)}{2m+1}}}\right),\num
\end{equation}
where $\psi$ is the $C^{\infty}(\rr,[0,1])$ function defined in (\ref{giu13}), and $w_2$ is a $C^{\infty}(\rr,[0,1])$ function satisfying similar properties as the function defined in (\ref{giu15}), with possibly different positive numerical values for its support localization, in order to have that 
\begin{equation}\label{bell1}\inc
\textrm{supp }\psi' \subset \big\{w_2=1\big\} \textrm{ and } \textrm{supp }w_2' \subset \big\{\psi=1\big\}. \num
\end{equation} 
In order to check (\ref{giu38}), we begin by noticing from (\ref{giu15}), (\ref{giu40}) and (\ref{giu41}) that for $0 \leq j \leq m-1$, 
\inc\begin{equation}\label{bell2}
r_{m-j-1}(X)^{\frac{1}{2m-2j-1}} \gtrsim r_{m-j}(X)^{\frac{1}{2m-2j+1}} \gtrsim ... \gtrsim 
r_{m-1}(X)^{\frac{1}{2m-1}} \gtrsim \langle X \rangle^{\frac{2}{2m+1}}, \num 
\end{equation}
on the support of the function
$$\textrm{supp} \Big(\tilde{W}_0\prod_{l=1}^{j}W_l\Big), \textrm{ if } 1 \leq j \leq m-1, \textrm{ or, } \textrm{supp }\tilde{W}_0,\textrm{ if } j=0.$$ 
Notice that the constants in the estimates (\ref{bell2}) only depend on the values of the parameters $\Lambda_0$,...,$\Lambda_{j-1}$ but not on $\Lambda_l$, when $l \geq j$.
This shows that the functions 
$$\Psi_0; \ \Big(\prod_{l=1}^{j}W_l\Big)\Psi_j, \textrm{ for } 1 \leq j \leq m-2; \textrm{ and }  \prod_{l=1}^{m-1}W_l,$$
are well-defined on the support of the function $\tilde{W}_0$. Now, by noticing from (\ref{giu13}), (\ref{giu15}), (\ref{giu39}), (\ref{giu40}) and (\ref{bell1}) that
\begin{equation}\label{eg1}\inc
1 \leq \Psi_j+W_{j+1},\num
\end{equation} 
on the support of the function
$$\textrm{supp} \Big(\tilde{W}_0\prod_{l=1}^{j}W_l\Big), \textrm{ if } 1 \leq j \leq m-2, \textrm{ or, } \textrm{supp }\tilde{W}_0,\textrm{ if } j=0,$$ 
we deduce the estimate (\ref{giu38}) from a finite iteration by using the following estimates
$$\tilde{W}_0 \leq \tilde{W}_0\Psi_0+\tilde{W}_0W_1$$
and
$$\tilde{W}_{0}\Big(\prod_{l=1}^{j}W_l\Big) \leq \tilde{W}_{0}\Big(\prod_{l=1}^{j}W_l\Big)\Psi_j+\tilde{W}_{0}\Big(\prod_{l=1}^{j+1}W_l\Big),$$
for any $1 \leq j \leq m-2$. One can also notice that (\ref{eg1}) implies that 
\begin{equation}\label{eg2}\inc
1 \leq \Psi_j+\sum_{k=j+1}^{m-2}{\Big(\prod_{l=j+1}^kW_l\Big)\Psi_k}+\prod_{l=j+1}^{m-1}W_l,\num
\end{equation}
on the support of the function
$$\textrm{supp} \Big(\tilde{W}_0\prod_{l=1}^{j}W_l\Big), \textrm{ if } 1 \leq j \leq m-2, \textrm{ or, } \textrm{supp }\tilde{W}_0,\textrm{ if } j=0.$$ 
Since $\textrm{Re }q_p \geq 0$, we then get from (\ref{bell2}) that 
\begin{equation}\label{giu42}\inc 
\forall X \in \rr^{2n}, \ \tilde{W}_{0}(X)\Big(\prod_{l=1}^{m-1}W_l(X)\Big)\langle X \rangle^{\frac{2}{2m+1}} \leq \tilde{a}_{\Lambda_0,...,\Lambda_{m-2}}  \sum_{p=1}^N\textrm{Re }q_p(X), \num
\end{equation}
where $\tilde{a}_{\Lambda_0,...,\Lambda_{m-2}}$ is a positive constant whose value depends on the parameters 
$$(\Lambda_l)_{0 \leq l \leq m-2}.$$
We define for $1 \leq p \leq N$,
\begin{equation}\label{giu44}\inc
\mathfrak{p}_{j,p}(X)=\tilde{W}_0(X)\Big(\prod_{l=1}^{j}{W_l(X)}\Big)\Psi_j(X)\frac{\tilde{r}_{m-j-1,p}(X)}{r_{m-j-1}(X)^{\frac{2m-2j-2}{2m-2j-1}}},  \num
\end{equation}
for $1 \leq j \leq m-2$, and
\begin{equation}\label{giu45}\inc
\mathfrak{p}_{0,p}(X)=\tilde{W}_0(X)\Psi_0(X)\frac{\tilde{r}_{m-1,p}(X)}{r_{m-1}(X)^{\frac{2m-2}{2m-1}}}, \num
\end{equation}
where the quadratic forms $\tilde{r}_{k,p}$ are defined in (\ref{giu11}) and (\ref{giu11serena}). We get from (\ref{giu13}), (\ref{giu15}), (\ref{giu39}), (\ref{giu40}), (\ref{giu41}), (\ref{bell2}), Lemma~\ref{lem3},  Lemma~\ref{lem2.1}, Lemma~\ref{lem2.2} and Lemma~\ref{lem2.3} that 
\begin{equation}\label{giu47}\inc
\mathfrak{p}_{j,p} \in S\big(1,\langle X \rangle^{-\frac{2(2m-2j-3)}{2m+1}}dX^2\big). \num
\end{equation}
for any $0 \leq j \leq m-2$.

We shall now study the Poisson brackets $H_{\textrm{Im}q_p} \mathfrak{p}_{j,p}$. In doing so, we begin by writing that
\inc\begin{align*}\label{giu48}
& \ H_{\textrm{Im}q_p}  \mathfrak{p}_{j,p}(X)=\big(H_{\textrm{Im}q_p}\tilde{W}_0\big)(X)\Big(\prod_{l=1}^{j}{W_l(X)}\Big)\Psi_j(X) \frac{\tilde{r}_{m-j-1,p}(X)}{r_{m-j-1}(X)^{\frac{2m-2j-2}{2m-2j-1}}} \num \\
+ & \ \tilde{W}_0(X)\Big(\prod_{l=1}^{j}{W_l(X)}\Big)\big(H_{\textrm{Im}q_p}\Psi_j\big)(X) \frac{\tilde{r}_{m-j-1,p}(X)}{r_{m-j-1}(X)^{\frac{2m-2j-2}{2m-2j-1}}} \\
+ & \  \tilde{W}_0(X)\Big(\prod_{l=1}^{j}{W_l(X)}\Big)\Psi_j(X)H_{\textrm{Im}q_p}\big(r_{m-j-1}(X)^{-\frac{2m-2j-2}{2m-2j-1}}\big)\tilde{r}_{m-j-1,p}(X) \\
+ & \ \tilde{W}_0(X)\Big(\prod_{l=1}^{j}{W_l(X)}\Big)\Psi_j(X)\frac{H_{\textrm{Im}q_p} \tilde{r}_{m-j-1,p}(X)}{r_{m-j-1}(X)^{\frac{2m-2j-2}{2m-2j-1}}}\\
+ & \ \sum_{l=1}^j\tilde{W}_0(X)\big(H_{\textrm{Im}q_p}W_l\big)(X)\Big(\prod_{k=1 \atop k \neq l}^{j}{W_k(X)}\Big)\Psi_j(X) \frac{\tilde{r}_{m-j-1,p}(X)}
{r_{m-j-1}(X)^{\frac{2m-2j-2}{2m-2j-1}}},
\end{align*}
for $1 \leq j \leq m-2$. We denote by respectively $B_{1,j,p}$, $B_{2,j,p}$, $B_{3,j,p}$, $B_{4,j,p}$ and $B_{5,j,p}$ the five terms appearing in the right hand side of (\ref{giu48}). 
We also write in the case where $j=0$,
\inc\begin{align*}\label{giu48bis}
& \ H_{\textrm{Im}q_p} \mathfrak{p}_{0,p}(X)=\big(H_{\textrm{Im}q_p}\tilde{W}_0\big)(X)\Psi_0(X) \frac{\tilde{r}_{m-1,p}(X)}{r_{m-1}(X)^{\frac{2m-2}{2m-1}}} \num \\
+ & \ \tilde{W}_0(X)\big(H_{\textrm{Im}q_p}\Psi_0\big)(X) \frac{\tilde{r}_{m-1,p}(X)}{r_{m-1}(X)^{\frac{2m-2}{2m-1}}} \\
+ & \  \tilde{W}_0(X)\Psi_0(X)H_{\textrm{Im}q_p}\big(r_{m-1}(X)^{-\frac{2m-2}{2m-1}}\big)\tilde{r}_{m-1,p}(X) \\
+ & \ \tilde{W}_0(X)\Psi_0(X)\frac{H_{\textrm{Im}q_p}\tilde{r}_{m-1,p}(X)}{r_{m-1}(X)^{\frac{2m-2}{2m-1}}},
\end{align*}
and denote as before by respectively $B_{1,0,p}$, $B_{2,0,p}$, $B_{3,0,p}$ and $B_{4,0,p}$ the four terms appearing in the right hand side of (\ref{giu48bis}).

Since the constants in the estimates (\ref{bell2}) only depend on the values of the parameters $\Lambda_0$,..., $\Lambda_{j-1}$; but not on $\Lambda_l$, when $l \geq j$; we notice from (\ref{giu38}), (\ref{bell2}) and (\ref{giu42}) that there exist $a_{0}>0$ and some positive constants $a_{j,\Lambda_0,...,\Lambda_{j-1}}$, for $1 \leq j \leq m-1$, whose values with respect to the parameters $(\Lambda_l)_{0 \leq l \leq m-2}$ only depend on $\Lambda_0$,..., $\Lambda_{j-1}$; but not on $\Lambda_l$, when $l \geq j$; such that for any constants $(\alpha_j)_{1 \leq j \leq m-2}$, with $\alpha_j \geq 1$; and $X \in \rr^{2n}$,
\inc\begin{align*}\label{giu100}
 \num & \ w_1\left(\frac{r_{m-1}(X)}{\langle X \rangle^{\frac{2(2m-1)}{2m+1}}}\right) \langle X \rangle^{\frac{2}{2m+1}} 
\leq a_0 \tilde{W}_0(X)\Psi_0(X) r_{m-1}(X)^{\frac{1}{2m-1}}\\ 
 & \  + \sum_{j=1}^{m-2}{\alpha_j a_{j,\Lambda_0,...,\Lambda_{j-1}} \tilde{W}_0(X)
\Big(\prod_{l=1}^{j}W_l(X)\Big)\Psi_j(X)}r_{m-j-1}(X)^{\frac{1}{2m-2j-1}}\\
  & \ +a_{m-1,\Lambda_0,...,\Lambda_{m-2}}\sum_{p=1}^N\textrm{Re }q_p(X).
\end{align*}
The positive constant $a_0$ is independent of any of the parameters $(\Lambda_l)_{0 \leq l \leq m-2}$.
Setting 
\begin{equation}\label{giu101}\inc
\mathfrak{p}_p=a_0 \mathfrak{p}_{0,p}+\sum_{j=1}^{m-2}{\alpha_j a_{j,\Lambda_0,...,\Lambda_{j-1}} \mathfrak{p}_{j,p}}, \num
\end{equation}
we know from (\ref{giu47}) that 
\begin{equation}\label{giu102}\inc
\mathfrak{p}_p \in S\big(1,\langle X \rangle^{-\frac{2}{2m+1}}dX^2\big).\num
\end{equation}
For any $\eps >0$, we shall prove that after a proper choice for the constants $(\Lambda_j)_{0 \leq j \leq m-2}$ and $(\alpha_j)_{1 \leq j \leq m-2}$, with $\Lambda_j \geq1$, $\alpha_j \geq 1$, whose values will depend on $\eps$; one can find a positive constant $c_{12,\eps}>0$ such that for all $X \in \rr^{2n}$,
\inc\begin{equation}\label{giu103}
 c_{12,\eps}\sum_{p=1}^N\big(\textrm{Re }q_p(X)+H_{\textrm{Im}q_p} \mathfrak{p}_p(X)\big) +\eps \langle X \rangle^{\frac{2}{2m+1}} \geq  w_1\left(\frac{r_{m-1}(X)}{\langle X \rangle^{\frac{2(2m-1)}{2m+1}}}\right) \langle X \rangle^{\frac{2}{2m+1}}.\num
\end{equation}
Once this estimate proved, Proposition~\ref{prop1} will directly follow from (\ref{giu33}), (\ref{giu35}), (\ref{giu102}) and (\ref{giu103}), if we choose the positive parameter $\eps$ sufficiently small and consider the weight functions
$$g_p=c_{13,\eps}G_{m,p}+c_{14,\eps}\mathfrak{p}_p, \ 1 \leq p \leq N,$$ 
after a suitable choice for the positive constants $c_{13,\eps}$ and $c_{14,\eps}$.

Let $\eps>0$, it therefore remains to choose properly these constants $(\Lambda_j)_{0 \leq j \leq m-2}$ and $(\alpha_j)_{1 \leq j \leq m-2}$, with $\Lambda_j \geq1$, $\alpha_j \geq 1$, in order to satisfy (\ref{giu103}).

Recalling from (\ref{kee1}) that for all $1 \leq p \leq N$ and $0 \leq s \leq m-2$,
\inc\begin{multline*}\label{giu104}
H_{\textrm{Im}q_p}\tilde{r}_{m-s-1,p}(X)=
2 \hspace{-1cm} \sum_{\substack{j=1,...,N \\ (l_1,...,l_{m-s-2}) \in \{1,...,N\}^{m-s-2}}}\hspace{-1cm} \textrm{Re }q_j(\textrm{Im }F_{l_1}...\textrm{Im }F_{l_{m-s-2}}\textrm{Im }F_{p}X)
\\ +  2\hspace{-1cm} \sum_{\substack{j=1,...,N \\ (l_1,...,l_{m-s-2}) \in \{1,...,N\}^{m-s-2}}} \hspace{-1cm} \textrm{Re }q_j(\textrm{Im }F_{l_1}...\textrm{Im }F_{l_{m-s-2}}X;\textrm{Im }F_{l_1}...\textrm{Im }F_{l_{m-s-2}}(\textrm{Im }F_{p})^2X),\num
\end{multline*}
one can notice by expanding the term 
$$2a_{m-1,\Lambda_0,...,\Lambda_{m-2}}\sum_{p=1}^N\textrm{Re }q_p+\sum_{p=1}^NH_{\textrm{Im}q_p}\mathfrak{p}_p,$$ 
by using (\ref{giu48}), (\ref{giu48bis}) and (\ref{giu101}) that the terms in
$$a_0 \sum_{p=1}^NB_{4,0,p}+\sum_{j=1}^{m-2}{\alpha_j a_{j,\Lambda_0,...,\Lambda_{j-1}} \left(\sum_{p=1}^NB_{4,j,p}\right)},$$
produced by the terms associated to 
$$\hspace{-1cm} \sum_{\substack{j=1,...,N \\ (l_1,...,l_{m-s-2}) \in \{1,...,N\}^{m-s-2}}}\hspace{-1cm} \textrm{Re }q_j(\textrm{Im }F_{l_1}...\textrm{Im }F_{l_{m-s-2}}\textrm{Im }F_{p}X),$$ 
while using (\ref{giu104}), give exactly two times the term 
\inc\begin{align*}\label{giu105}
& \ a_0\tilde{W}_0(X)\Psi_0(X) r_{m-1}(X)^{\frac{1}{2m-1}} \num \\ 
+ & \  \sum_{j=1}^{m-2}{\alpha_j a_{j,\Lambda_0,...,\Lambda_{j-1}} \tilde{W}_0(X)
\Big(\prod_{l=1}^{j}W_l(X)\Big)\Psi_j(X)}r_{m-j-1}(X)^{\frac{1}{2m-2j-1}}\\
+ & \ a_{m-1,\Lambda_0,...,\Lambda_{m-2}}\sum_{p=1}^N\textrm{Re }q_p(X),
\end{align*}
for which we have the estimate (\ref{giu100}). To prove the estimate (\ref{giu103}), it will therefore be sufficient to check that all the other terms appearing in (\ref{giu48}) and (\ref{giu48bis}) can also be all absorbed in the term (\ref{giu105}) after a proper choice for the constants $(\Lambda_j)_{0 \leq j \leq m-2}$ and $(\alpha_j)_{1 \leq j \leq m-2}$; at the exception of a remainder term in 
$$\eps \langle X \rangle^{\frac{2}{2m+1}}.$$ 
We shall choose these constants in the following order $\Lambda_0$, $\alpha_1$, $\Lambda_1$, $\alpha_2$, ...., $\alpha_{m-2}$ and $\Lambda_{m-2}$. 

We successively study the remaining terms in (\ref{giu48}) and (\ref{giu48bis}), by increasing value of the integer $0 \leq j \leq m-2$. We first notice from (\ref{giu13}), (\ref{giu15}), (\ref{giu39}), (\ref{giu41}), (\ref{giu48bis}), Lemma~\ref{lem2.331} and Lemma~\ref{lem2.435} that one can choose the first constant $\Lambda_0 \geq 1$ such that for all $X \in \rr^{2n}$,
\begin{equation}\label{giu106}\inc
a_0 \sum_{p=1}^N|B_{1,0,p}(X)| \lesssim \Lambda_0^{-\frac{1}{2}} \langle X \rangle^{\frac{2}{2m+1}} \leq \frac{\eps}{m-1} \langle X \rangle^{\frac{2}{2m+1}}.\num
\end{equation}

By noticing from (\ref{bell2}) that the estimates
\begin{equation}\label{giu107}\inc
r_{m}(X) \lesssim \langle X \rangle^2 \lesssim r_{m-1}(X)^{\frac{2m+1}{2m-1}},\num
\end{equation}
are fulfilled on the support of the function $\tilde{W}_0$, we deduce from (\ref{giu13}), (\ref{giu39}) and (\ref{giu48bis}) that the modulus of the terms $B_{3,0,p}$ can be estimated as
\begin{align*}
a_{0}\sum_{p=1}^N|B_{3,0,p}(X)| 
= & \     a_{0}\sum_{p=1}^N \big| r_{m-1}(X)^{\frac{2m-2}{2m-1}} H_{\textrm{Im}q_p}\big(r_{m-1}(X)^{-\frac{2m-2}{2m-1}}\big)\big| \\
& \quad \quad \quad \times \big|r_{m-1}(X)^{-\frac{2m-2}{2m-1}} \tilde{r}_{m-1,p}(X)\big| \tilde{W}_0(X)\Psi_0(X) \\ 
\lesssim & \ \Lambda_0^{-\frac{1}{2}}\tilde{W}_0(X)\Psi_0(X)r_{m-1}(X)^{\frac{1}{2m-1}},
\end{align*}
for all $X \in \rr^{2n}$; 
since from Lemma~\ref{lem2.331} and Lemma~\ref{lem2.41}, we have for any $p$ in $\{1,...,N\}$ that 
$$ \big|r_{m-1}(X)^{\frac{2m-2}{2m-1}} H_{\textrm{Im}q_p}\big(r_{m-1}(X)^{-\frac{2m-2}{2m-1}}\big)\big|  \lesssim r_{m-1}(X)^{\frac{1}{2m-1}}$$
and 
$$\big|r_{m-1}(X)^{-\frac{2m-2}{2m-1}} \tilde{r}_{m-1,p}(X)\big| \lesssim \Lambda_0^{-\frac{1}{2}},$$
on the support of the function $\tilde{W}_0(X)\Psi_0(X).$
By possibly increasing sufficiently the value of the constant $\Lambda_0$ which is of course possible while keeping (\ref{giu106}), one can control this term with the \og good \fg \ term (\ref{giu105}).

Next, we deduce from (\ref{giu13}), (\ref{giu39}), (\ref{giu48bis}), (\ref{giu107}) and Lemma~\ref{lem2.4} that the modulus of the second terms in $B_{4,0,p}$ associated to
$$2\hspace{-1cm} \sum_{\substack{j=1,...,N \\ (l_1,...,l_{m-2}) \in \{1,...,N\}^{m-2}}} \hspace{-1cm} \textrm{Re }q_j(\textrm{Im }F_{l_1}...\textrm{Im }F_{l_{m-2}}X;\textrm{Im }F_{l_1}...\textrm{Im }F_{l_{m-2}}(\textrm{Im }F_{p})^2X),$$ 
while using (\ref{giu104}), denoted here $\tilde{B}_{4,0,p}$, 
\begin{align*}
\sum_{p=1}^N& \  \tilde{B}_{4,0,p}(X) =\tilde{W}_0(X)\Psi_0(X) \\
& \ \times  \sum_{p=1}^N\left(\frac{H_{\textrm{Im}q_p}\tilde{r}_{m-1,p}(X)}{r_{m-1}(X)^{\frac{2m-2}{2m-1}}}-
 2 \hspace{-1cm} \sum_{\substack{j=1,...,N \\ (l_1,...,l_{m-2}) \in \{1,...,N\}^{m-2}}}\hspace{-1cm} 
 \frac{\textrm{Re }q_j(\textrm{Im }F_{l_1}...\textrm{Im }F_{l_{m-2}}\textrm{Im }F_{p}X)}{r_{m-1}(X)^{\frac{2m-2}{2m-1}}}\right),\\
& \  = \tilde{W}_0(X)\Psi_0(X) \left(\sum_{p=1}^N\frac{H_{\textrm{Im}q_p}\tilde{r}_{m-1,p}(X)}{r_{m-1}(X)^{\frac{2m-2}{2m-1}}}-
 2 r_{m-1}(X)^{\frac{1}{2m-1}}\right)
\end{align*}
can be estimated as
$$a_{0}\sum_{p=1}^N|\tilde{B}_{4,0,p}(X)| 
\lesssim \Lambda_0^{-\frac{1}{2}}\tilde{W}_0(X)\Psi_0(X)r_{m-1}(X)^{\frac{1}{2m-1}},$$
for all $X \in \rr^{2n}$. By possibly increasing sufficiently the value of the constant $\Lambda_0$ which is of course possible while keeping (\ref{giu106}), one can also control this term with the \og good \fg \ term (\ref{giu105}). The value of the constant $\Lambda_0$ is now definitively fixed. In (\ref{giu48bis}), it only remains to study the terms $B_{2,0,p}$.

About these terms, we deduce from (\ref{giu13}), (\ref{giu39}), (\ref{giu48bis}), (\ref{giu107}), Lemma \ref{lem2.331} and Lemma~\ref{lem2.42} that for all $X \in \rr^{2n}$,
\begin{equation}\label{eg3}\inc
a_0\sum_{p=1}^N |B_{2,0,p}(X)| \lesssim \tilde{W}_0(X)W_1(X)r_{m-1}(X)^{\frac{1}{2m-1}}.\num
\end{equation}
By using now (\ref{bell2}) and (\ref{eg2}) with $j=1$, we obtain that for all $X \in \rr^{2n}$,
\begin{multline*}
a_0 \sum_{p=1}^N |B_{2,0,p}(X)| \leq c_{m-1,\Lambda_0,...,\Lambda_{m-2}}\tilde{W}_0(X)\Big(\prod_{l=1}^{m-1}W_l(X)\Big)\sum_{p=1}^N \textrm{Re }q_p(X) \\
+ \sum_{j=1}^{m-2}{c_{j,\Lambda_0,...,\Lambda_{j-1}}\tilde{W}_0(X)\Big(\prod_{l=1}^jW_l(X)\Big)\Psi_j(X)r_{m-j-1}(X)^{\frac{1}{2m-2j-1}}},
\end{multline*}
which implies that 
\inc\begin{multline*}\label{eg4}
a_0\sum_{p=1}^N |B_{2,0,p}(X)| \leq c_{m-1,\Lambda_0,...,\Lambda_{m-2}}\sum_{p=1}^N \textrm{Re }q_p(X)\\
+ \sum_{j=1}^{m-2}{c_{j,\Lambda_0,...,\Lambda_{j-1}}\tilde{W}_0(X)\Big(\prod_{l=1}^jW_l(X)\Big)\Psi_j(X)r_{m-j-1}(X)^{\frac{1}{2m-2j-1}}},\num
\end{multline*}
where the quantities $c_{j,\Lambda_0,...,\Lambda_{j-1}}$ stand for positive constants whose values depend on $\Lambda_0$,..., $\Lambda_{j-1}$, but not on 
$(\Lambda_k)_{j \leq k \leq m-2}$ and $(\alpha_k)_{1 \leq k \leq m-2},$ 
according to the remark done after (\ref{bell2}). One can therefore choose the constant $\alpha_1 \geq 1$ in (\ref{giu101}) sufficiently large in order to absorb the term of the index $j=1$ in the sum appearing in the right hand side of the estimate (\ref{eg4}) by the term of same index in the \og good\fg \ term (\ref{giu105}). This is possible since the constants $a_{1,\Lambda_0}$ and 
$c_{1,\Lambda_0}$ are now fixed after our choice of the parameter $\Lambda_0$.

This ends our step index $j=0$ in which we have chosen the values for the two constants $\Lambda_0$ and $\alpha_1 \geq 1$. We shall now explain how to choose the remaining constants
$(\Lambda_j)_{1 \leq j \leq m-2}$ and $(\alpha_j)_{2 \leq j \leq m-2}$ in (\ref{giu101}) in order to satisfy (\ref{giu103}). This choice will also determine the values of the constants $(a_{j,\Lambda_0,...,\Lambda_{j-1}})_{1 \leq j \leq m-2}$ appearing in (\ref{giu101}). After this step index $j=0$, we have managed to absorb all the terms appearing in (\ref{giu48bis}) in the \og good \fg \ term (\ref{giu105}) at the exception of a remainder coming from (\ref{giu106}) and (\ref{eg4}),
$$\sum_{j=2}^{m-2}{c_{j,\Lambda_0,...,\Lambda_{j-1}}\tilde{W}_0(X)\Big(\prod_{l=1}^jW_l(X)\Big)\Psi_j(X)r_{m-j-1}(X)^{\frac{1}{2m-2j-1}}}
+\frac{\eps}{m-1}\langle X\rangle^{\frac{2}{2m+1}},$$
where one recall that the positive constants $c_{j,\Lambda_0,...,\Lambda_{j-1}}$ only depend on $\Lambda_0$,...,$\Lambda_{j-1}$, but not on 
$(\Lambda_k)_{j \leq k \leq m-2}$ and $(\alpha_k)_{1 \leq k \leq m-2}.$

We proceed in the following by finite induction and assume that, at the beginning of the step index $k$, with $1 \leq k \leq m-2$, we have already chosen the values for the constants 
$(\Lambda_j)_{0 \leq j \leq k-1}$ and $(\alpha_j)_{1 \leq j \leq k}$ in (\ref{giu101}); and that these choices have allowed to absorb all the terms appearing in the right hand side of (\ref{giu48bis}) and (\ref{giu48}), when $1 \leq j \leq k-1$, in the \og good \fg \ term (\ref{giu105}) at the exception of a remainder term
\inc\begin{multline*}\label{eg5}
\frac{k}{m-1}\eps \langle X \rangle^{\frac{2}{2m+1}}+\\
\sum_{j=k+1}^{m-2}\tilde{c}_{j,\Lambda_0,...,\Lambda_{j-1}, \alpha_1,...,\alpha_{k-1}}\tilde{W}_0(X)\Big(\prod_{l=1}^jW_l(X)\Big)\Psi_j(X)r_{m-j-1}(X)^{\frac{1}{2m-2j-1}},\num
\end{multline*} 
where the quantities $\tilde{c}_{j,\Lambda_0,...,\Lambda_{j-1}, \alpha_1,...,\alpha_{k-1}}$ stand for positive constants whose values only depend on $\Lambda_0$,..., $\Lambda_{j-1}$, $\alpha_1$,..., $\alpha_{k-1}$; but not on $(\Lambda_l)_{j \leq l \leq m-2}$ and $(\alpha_l)_{k \leq l \leq m-2}$.

We shall now explain how to choose the constants $\Lambda_k$ and; $\alpha_{k+1}$, when $k \leq m-3$; in this step index $k$ in order to absorb the terms appearing in the right hand side of (\ref{giu48}), when $j=k$, at the exception of a remainder term of the type (\ref{eg5}) where $k$ will be replaced by $k+1$;  in the \og good \fg \ term (\ref{giu105}). Since the constants $(\Lambda_j)_{0 \leq j \leq k-1}$ and $(\alpha_j)_{1 \leq j \leq k}$ have already been chosen, we shall only underline in the following the dependence of our estimates with respect to the other parameters 
$(\Lambda_j)_{k \leq j \leq m-2}$ and $(\alpha_j)_{k+1 \leq j \leq m-2}$, whose values remain to be chosen.

We notice from (\ref{giu13}), (\ref{giu39}), (\ref{giu40}), (\ref{giu41}), (\ref{bell2}), (\ref{giu48}), Lemma~\ref{lem2.331} and Lemma \ref{lem2.435} that one can assume by choosing the constant $\Lambda_k \geq 1$ sufficiently large that for all $X \in \rr^{2n}$,
\begin{equation}\label{eg6}\inc
\alpha_k a_{k,\Lambda_0,...,\Lambda_{k-1}} \sum_{p=1}^N |B_{1,k,p}(X)| \lesssim \Lambda_k^{-\frac{1}{2}} \langle X \rangle^{\frac{2}{2m+1}} \leq \frac{\eps}{m-1} \langle X \rangle^{\frac{2}{2m+1}},\num
\end{equation}
since the constants $\alpha_k$, $\Lambda_0$,....,$\Lambda_{k-1}$ have already been fixed.

Next, we deduce from (\ref{giu13}), (\ref{giu39}), (\ref{bell2}) and (\ref{giu48}) that the modulus of the terms $B_{3,k,p}$ can be estimated as
\begin{align*}
& \ \alpha_k a_{k,\Lambda_0,...,\Lambda_{k-1}}\sum_{p=1}^N|B_{3,k,p}(X)| \\
= & \     \alpha_k a_{k,\Lambda_0,...,\Lambda_{k-1}}\sum_{p=1}^N \big| r_{m-k-1}(X)^{\frac{2m-2k-2}{2m-2k-1}} H_{\textrm{Im}q_p}\big(r_{m-k-1}(X)^{-\frac{2m-2k-2}{2m-2k-1}}\big)\big| \\
& \quad \quad \quad \times \big|r_{m-k-1}(X)^{-\frac{2m-2k-2}{2m-2k-1}} \tilde{r}_{m-k-1,p}(X)\big| \tilde{W}_0(X)\Big(\prod_{l=1}^kW_l(X)\Big)\Psi_k(X) \\ 
\lesssim & \ \Lambda_k^{-\frac{1}{2}}\tilde{W}_0(X)\Big(\prod_{l=1}^kW_l(X)\Big)\Psi_k(X)r_{m-k-1}(X)^{\frac{1}{2m-2k-1}},
\end{align*}
for all $X \in \rr^{2n}$; 
since from Lemma~\ref{lem2.331} and Lemma~\ref{lem2.41}, we have for any $p$ in $\{1,...,N\}$ that 
$$ \big|r_{m-k-1}(X)^{\frac{2m-2k-2}{2m-2k-1}} H_{\textrm{Im}q_p}\big(r_{m-k-1}(X)^{-\frac{2m-2k-2}{2m-2k-1}}\big)\big|  \lesssim r_{m-k-1}(X)^{\frac{1}{2m-2k-1}}$$
and 
$$\big|r_{m-k-1}(X)^{-\frac{2m-2k-2}{2m-2k-1}} \tilde{r}_{m-k-1,p}(X)\big| \lesssim \Lambda_k^{-\frac{1}{2}},$$
on the support of the function 
$$\tilde{W}_0(X)\Big(\prod_{l=1}^kW_l(X)\Big)\Psi_k(X).$$
By possibly increasing sufficiently the value of the constant $\Lambda_k$ which is of course possible while keeping (\ref{eg6}), one can control this term with the \og good \fg \ term (\ref{giu105}).

Next, we deduce from (\ref{giu13}), (\ref{giu39}), (\ref{bell2}), (\ref{giu48}) and Lemma~\ref{lem2.4} that the modulus of the second terms in $B_{4,k,p}$ associated to
$$2\hspace{-1cm} \sum_{\substack{j=1,...,N \\ (l_1,...,l_{m-k-2}) \in \{1,...,N\}^{m-k-2}}} \hspace{-1cm} \textrm{Re }q_j(\textrm{Im }F_{l_1}...\textrm{Im }F_{l_{m-k-2}}X;\textrm{Im }F_{l_1}...\textrm{Im }F_{l_{m-k-2}}(\textrm{Im }F_{p})^2X),$$ 
while using (\ref{giu104}), denoted here $\tilde{B}_{4,k,p}$, 
\begin{align*}
\sum_{p=1}^N& \  \tilde{B}_{4,k,p}(X) =\tilde{W}_0(X)\Big(\prod_{l=1}^kW_l(X)\Big)\Psi_k(X) \\
& \ \times  \sum_{p=1}^N\left(\frac{H_{\textrm{Im}q_p}\tilde{r}_{m-k-1,p}(X)}{r_{m-k-1}(X)^{\frac{2m-2k-2}{2m-2k-1}}}-
 2 \hspace{-1cm} \sum_{\substack{j=1,...,N \\ (l_1,...,l_{m-k-2}) \in \{1,...,N\}^{m-k-2}}}\hspace{-1cm} 
 \frac{\textrm{Re }q_j(\textrm{Im }F_{l_1}...\textrm{Im }F_{l_{m-k-2}}\textrm{Im }F_{p}X)}{r_{m-k-1}(X)^{\frac{2m-2k-2}{2m-2k-1}}}\right),\\
& \  = \tilde{W}_0(X)\Big(\prod_{l=1}^kW_l(X)\Big)\Psi_k(X) \left(\sum_{p=1}^N\frac{H_{\textrm{Im}q_p}\tilde{r}_{m-k-1,p}(X)}{r_{m-k-1}(X)^{\frac{2m-2k-2}{2m-2k-1}}}-
 2 r_{m-k-1}(X)^{\frac{1}{2m-2k-1}}\right)
\end{align*}
can be estimated as
$$\alpha_k a_{k,\Lambda_0,...,\Lambda_{k-1}}\sum_{p=1}^N|\tilde{B}_{4,k,p}(X)| 
\lesssim \Lambda_k^{-\frac{1}{2}}\tilde{W}_0(X)\Big(\prod_{l=1}^kW_l(X)\Big)\Psi_k(X)r_{m-k-1}(X)^{\frac{1}{2m-2k-1}},$$
for all $X \in \rr^{2n}$. By possibly increasing sufficiently the value of the constant $\Lambda_k$ which is of course possible while keeping (\ref{eg6}), one can also control this term with the \og good \fg \ term (\ref{giu105}).

For $1 \leq l \leq k$ and $1 \leq p \leq N$, we shall now study the term
$$B_{5,k,p,l}(X)=\tilde{W}_0(X)\big(H_{\textrm{Im}q_p}W_l\big)(X)\Big(\prod_{j=1 \atop j \neq l}^{k}{W_j(X)}\Big)\Psi_k(X) \frac{\tilde{r}_{m-k-1,p}(X)}{r_{m-k-1}(X)^{\frac{2m-2k-2}{2m-2k-1}}},$$
appearing in the term $B_{5,k,p}$ in (\ref{giu48}). 
By noticing that 
$$r_{m-l-2}(X) \sim \Lambda_l^{-1}r_{m-l-1}(X)^{\frac{2m-2l-3}{2m-2l-1}},$$
on the support of the function $H_{\textrm{Im}q_p}W_{l+1}$, it follows from (\ref{giu13}), (\ref{giu15}), (\ref{giu39}), (\ref{giu40}), (\ref{giu41}), (\ref{bell2}), (\ref{giu107}), Lemma \ref{lem2.331}
and Lemma~\ref{lem2.43} that for all $X \in \rr^{2n}$,
$$\alpha_k a_{k,\Lambda_0,...,\Lambda_{k-1}}\sum_{p=1}^N |B_{5,k,p,1}(X)| \lesssim \Lambda_k^{-\frac{1}{2}}\tilde{W}_0(X)\Psi_{0}(X)r_{m-1}(X)^{\frac{1}{2m-1}}$$
and
$$\alpha_k a_{k,\Lambda_0,...,\Lambda_{k-1}}\sum_{p=1}^N |B_{5,k,p,l}(X)|  
\lesssim \Lambda_k^{-\frac{1}{2}}\tilde{W}_0(X)\Big(\prod_{j=1}^{l-1}{W_j(X)}\Big)\Psi_{l-1}(X)r_{m-l}(X)^{\frac{1}{2m-2l+1}},$$
when $l \geq 2$. By possibly increasing again the value of the constant $\Lambda_k$, one can therefore control the term 
$$\alpha _k a_{k,\Lambda_0,...,\Lambda_{k-1}} \sum_{p=1}^N B_{5,k,p},$$
 with the \og good \fg \ term (\ref{giu105}).
The value of the constant $\Lambda_k$ is now definitively fixed.

About the terms $B_{2,k,p}$, we deduce from (\ref{giu13}), (\ref{giu39}), (\ref{bell2}), (\ref{giu48}), Lemma \ref{lem2.331} and Lemma~\ref{lem2.42} that for all $X \in \rr^{2n}$,
\inc\begin{equation}\label{eg7}
\alpha_k a_{k,\Lambda_0,...,\Lambda_{k-1}}\sum_{p=1}^N |B_{2,k,p}(X)|  \lesssim \tilde{W}_0(X)\Big(\prod_{l=1}^{k+1}W_l(X)\Big)r_{m-k-1}(X)^{\frac{1}{2m-2k-1}}.\num
\end{equation}
By distinguishing two cases, we first assume in the following that $k \leq m-3$. In this case,
by using (\ref{bell2}) and (\ref{eg2}) with $j=k+1$, we obtain that for all $X \in \rr^{2n}$,
\begin{align*}
& \ \alpha_k a_{k,\Lambda_0,...,\Lambda_{k-1}}\sum_{p=1}^N |B_{2,k,p}(X)| \\
\leq & \  c_{m-1,\Lambda_0,...,\Lambda_{m-2},\alpha_1,...,\alpha_k}'\tilde{W}_0(X)\Big(\prod_{l=1}^{m-1}W_l(X)\Big)\sum_{p=1}^N\textrm{Re }q_p(X)\\
& \ +\sum_{j=k+1}^{m-2}{c_{j,\Lambda_0,...,\Lambda_{j-1},\alpha_1,...,\alpha_k}'\tilde{W}_0(X)\Big(\prod_{l=1}^jW_l(X)\Big)\Psi_j(X)r_{m-j-1}(X)^{\frac{1}{2m-2j-1}}},
\end{align*}
which implies that 
\inc\begin{multline*}\label{eg8}
\alpha_k a_{k,\Lambda_0,...,\Lambda_{k-1}}\sum_{p=1}^N |B_{2,k,p}(X)| \leq  c_{m-1,\Lambda_0,...,\Lambda_{m-2},\alpha_1,...,\alpha_k}'\sum_{p=1}^N\textrm{Re }q_p(X)\\
+\sum_{j=k+1}^{m-2}{c_{j,\Lambda_0,...,\Lambda_{j-1},\alpha_1,...,\alpha_k}'\tilde{W}_0(X)\Big(\prod_{l=1}^jW_l(X)\Big)\Psi_j(X)r_{m-j-1}(X)^{\frac{1}{2m-2j-1}}},\num
\end{multline*}
where the quantities $c_{j,\Lambda_0,...,\Lambda_{j-1},\alpha_1,...,\alpha_k}'$ stand for positive constants whose values only depend on $\Lambda_0$,..., $\Lambda_{j-1}$, $\alpha_1$,..., $\alpha_{k}$, but not on $(\Lambda_l)_{j \leq l \leq m-2}$ and $(\alpha_l)_{k+1 \leq l \leq m-2}$. Indeed, we recall that the constants appearing in the estimates (\ref{bell2}) only depend on the values of the parameters $\Lambda_0$,..., $\Lambda_{j-1}$; but not on $(\Lambda_l)_{j \leq l \leq m-2}$ and $(\alpha_l)_{1 \leq l \leq m-2}$. One can therefore choose the constant $\alpha_{k+1} \geq 1$ in (\ref{giu101}) sufficiently large in order to absorb the term of index $j=k+1$ in the sum (\ref{eg5});
and the term of index $j=k+1$ in the sum appearing in the right hand side of the estimate (\ref{eg8}), by the term of same index in the \og good\fg \ term (\ref{giu105}).

When $k=m-2$ and taking $\Lambda_{m-2}=1$, it follows from (\ref{bell2}), used with $j=m-1$, and (\ref{eg7}) that for all $X \in \rr^{2n}$,
\inc\begin{align*}\label{eg9}
\alpha_{m-2} a_{m-2,\Lambda_0,...,\Lambda_{m-3}} \sum_{p=1}^N |B_{2,m-2,p}(X)| \lesssim & \ \tilde{W}_0(X)\Big(\prod_{l=1}^{m-1}W_l(X)\Big)r_1(X)^{\frac{1}{3}} \num \\
\lesssim & \  \sum_{p=1}^N \textrm{Re }q_p(X).
\end{align*}
This process allows us to achieve the construction of the weight functions $\mathfrak{p}_p$, $1 \leq p \leq N$, satisfying (\ref{giu103}), which ends the proof of (\ref{giu103}). This also ends the proof of Proposition~\ref{prop1}.~$\Box$

\section{Appendix}\label{appendix}
\init

\subsection{Wick calculus}\label{wick} The purpose of this section is to recall the definition and basic properties of the Wick quantization that we need for the proof of Theorem~\ref{theorem1}. We follow here the presentation of the Wick quantization given by N.~Lerner in \cite{lerner} and refer the reader to his work for the proofs of the results recalled below.

The main property of the Wick quantization is its property of positivity, i.e., that non-negative Hamiltonians define non-negative operators
$$a \geq 0 \Rightarrow a^{\textrm{Wick}} \geq 0.$$
We recall that this is not the case for the Weyl quantization and refer to \cite{lerner} for an explicit example of non-negative Hamiltonian defining an operator which is not non-negative.

Before defining properly the Wick quantization, we first need to recall the definition of the wave packets transform of a function $u \in \mathcal{S}(\rr^n)$, 
$$Wu(y,\eta)=(u,\varphi_{y,\eta})_{L^2(\rr^n)}=2^{n/4}\int_{\rr^n}{u(x)e^{- \pi (x-y)^2}e^{-2i \pi(x-y).\eta}dx}, \ (y,\eta) \in \rr^{2n}.$$
where 
$$\varphi_{y,\eta}(x)=2^{n/4}e^{- \pi (x-y)^2}e^{2i \pi (x-y).\eta}, \ x \in \mathbb{R}^n,$$
and $x^2=x_1^2+...+x_n^2$. With this definition, one can check (see Lemma 2.1 in \cite{lerner}) that 
the mapping $u \mapsto Wu$ is continuous from $\mathcal{S}(\rr^n)$ to $\mathcal{S}(\rr^{2n})$, isometric from $L^{2}(\rr^n)$ to $L^2(\rr^{2n})$ and that we have the
reconstruction formula
\begin{equation}\label{lay0.1}\inc
\forall u \in \mathcal{S}(\rr^n), \forall x \in \rr^n, \ u(x)=\int_{\rr^{2n}}{Wu(y,\eta)\varphi_{y,\eta}(x)dyd\eta}.\num
\end{equation}
By denoting $\Sigma_Y$ the operator defined in the Weyl quantization by the symbol 
$$p_Y(X)=2^n e^{-2\pi|X-Y|^2}, \ Y=(y,\eta) \in \rr^{2n},$$
which is a rank-one orthogonal projection
$$\big{(}\Sigma_Y u\big{)}(x)=Wu(Y)\varphi_Y(x)=(u,\varphi_Y)_{L^2(\rr^n)}\varphi_Y(x),$$
we define the Wick quantization of any $L^{\infty}(\rr^{2n})$  symbol $a$ as
\begin{equation}\label{lay0.2}\inc
a^{\textrm{Wick}}=\int_{\rr^{2n}}{a(Y)\Sigma_Y dY}.\num
\end{equation}
More generally, one can extend this definition when the symbol $a$ belongs to $\mathcal{S}'(\rr^{2n})$ by defining the operator $a^{\textrm{Wick}}$ for any $u$ and $v$ in $\mathcal{S}(\rr^{n})$ by
$$<a^{\textrm{Wick}}u,\overline{v}>_{\mathcal{S}'(\rr^{n}),\mathcal{S}(\rr^{n})}=<a(Y),(\Sigma_Yu,v)_{L^2(\rr^n)}>_{\mathcal{S}'(\rr^{2n}),\mathcal{S}(\rr^{2n})},$$
where $<\textrm{\textperiodcentered},\textrm{\textperiodcentered}>_{\mathcal{S}'(\rr^n),\mathcal{S}(\rr^n)}$ denotes the duality bracket between the
spaces $\mathcal{S}'(\rr^n)$ and $\mathcal{S}(\rr^n)$. The Wick quantization is a positive quantization
\begin{equation}\label{lay0.5}\inc
a \geq 0 \Rightarrow a^{\textrm{Wick}} \geq 0. \num
\end{equation}
In particular, real Hamiltonians get quantized in this quantization by formally self-adjoint operators and one has (see Proposition 3.2 in \cite{lerner}) that $L^{\infty}(\rr^{2n})$ symbols define bounded operators on $L^2(\rr^n)$ such that 
\begin{equation}\label{lay0}\inc
\|a^{\textrm{Wick}}\|_{\mathcal{L}(L^2(\rr^n))} \leq \|a\|_{L^{\infty}(\rr^{2n})}.\num
\end{equation}
According to Proposition~3.3 in~\cite{lerner}, the Wick and Weyl quantizations of a symbol $a$ are linked by the following identities
\begin{equation}\label{lay1bis}\inc
a^{\textrm{Wick}}=\tilde{a}^w,\num
\end{equation}
with
\begin{equation}\label{lay2bis}\inc
\tilde{a}(X)=\int_{\rr^{2n}}{a(X+Y)e^{-2\pi |Y|^2}2^ndY}, \ X \in \rr^{2n},\num
\end{equation}
and
\begin{equation}\label{lay1}\inc
a^{\textrm{Wick}}=a^w+r(a)^w,\num
\end{equation}
where $r(a)$ stands for the symbol
\begin{equation}\label{lay2}\inc
r(a)(X)=\int_0^1\int_{\rr^{2n}}{(1-\theta)a''(X+\theta Y)Y^2e^{-2\pi |Y|^2}2^ndYd\theta}, \ X \in \rr^{2n},\num
\end{equation}
if we use here the normalization chosen in \cite{lerner} for the Weyl quantization
\begin{equation}\label{lay3}\inc
(a^wu)(x)=\int_{\rr^{2n}}{e^{2i\pi(x-y).\xi}a\Big(\frac{x+y}{2},\xi\Big)u(y)dyd\xi},\num
\end{equation}
which differs from the one chosen in this paper. Because of this difference in normalizations, certain constant factors will naturally appear in the core of the proof of Theorem~\ref{theorem1}
while using certain formulas of Section~\ref{wick}, but these are minor adaptations.  
We also recall the following composition formula obtained in the proof of Proposition~3.4 in~\cite{lerner},
\begin{equation}\label{lay4}\inc
a^{\textrm{Wick}} b^{\textrm{Wick}} =\Big{[}ab-\frac{1}{4 \pi} a' \cdot b'+\frac{1}{4i \pi}\{a,b\} \Big{]}^{\textrm{Wick}}+S, \num
\end{equation}
with $\|S\|_{\mathcal{L}(L^2(\rr^n))} \leq d_n \|a\|_{L^{\infty}}\gamma_{2}(b),$
when $a \in L^{\infty}(\rr^{2n})$ and $b$ is a smooth symbol satisfying
$$\gamma_2(b)=\sup_{X \in \rr^{2n}, \atop T \in \rr^{2n}, |T|=1}|b^{(2)}(X)T^2| < +\infty.$$ 
The term $d_n$ appearing in the previous estimate stands for a positive constant depending only on the dimension $n$, and the notation $\{a,b\}$ denotes the Poisson bracket
$$\{a,b\}=\frac{\partial a}{\partial \xi} \cdot \frac{\partial b}{\partial x}-\frac{\partial a}{\partial x} \cdot \frac{\partial b}{\partial \xi}.$$

\subsection{Some technical lemmas}

This second part of the appendix is devoted to the proofs of several technical lemmas.

\bigskip

\begin{lemma}\label{lem2}
For any $1 \leq j \leq N$, $1 \leq p \leq N$, $(l_1,...,l_{k}) \in \{1,...,N\}^k$ and $s_1, s_2 \in \nn$, we have
\inc\begin{align*}\label{mari5}
& \ H_{\emph{\textrm{Im}}q_p} \Big( \emph{\textrm{Re }}q_j\big(\emph{\textrm{Im }}F_{l_1}...\emph{\textrm{Im }}F_{l_{k}}(\emph{\textrm{Im }}F_{p})^{s_1}X;\emph{\textrm{Im }}F_{l_1}...\emph{\textrm{Im }}F_{l_{k}}(\emph{\textrm{Im }}F_{p})^{s_2}X\big)\Big) \num \\
= &\ 2\emph{\textrm{Re }}q_j\big(\emph{\textrm{Im }}F_{l_1}...\emph{\textrm{Im }}F_{l_{k}}(\emph{\textrm{Im }}F_{p})^{s_1+1}X;\emph{\textrm{Im }}F_{l_1}...\emph{\textrm{Im }}F_{l_{k}}(\emph{\textrm{Im }}F_{p})^{s_2}X\big)\\ 
+ & \ 2 \emph{\textrm{Re }}q_j\big(\emph{\textrm{Im }}F_{l_1}...\emph{\textrm{Im }}F_{l_{k}}(\emph{\textrm{Im }}F_{p})^{s_1}X;\emph{\textrm{Im }}F_{l_1}...\emph{\textrm{Im }}F_{l_{k}}(\emph{\textrm{Im }}F_{p})^{s_2+1}X\big),
\end{align*}
where $\emph{\textrm{Re }}q_j(X;Y)$ stands for the polarized form associated to the quadratic form $\emph{\textrm{Re }}q_j$.
\end{lemma}

\bigskip

\noindent
\textit{Proof of Lemma~\ref{lem2}}. 
We begin by noticing from (\ref{10}) and the skew-symmetry property of Hamilton maps (\ref{12}) that the Hamilton map of the quadratic form
$$\tilde{r}(X)=\textrm{Re }q_j\big(\textrm{Im }F_{l_1}...\textrm{Im }F_{l_{k}}(\textrm{Im }F_{p})^{s_1}X;\textrm{Im }F_{l_1}...\textrm{Im }F_{l_{k}}(\textrm{Im }F_{p})^{s_2}X\big),$$ 
is given by 
\inc\begin{multline*}\label{lol10}
\tilde{F}=\frac{1}{2}(-1)^{k+s_1}(\textrm{Im }F_{p})^{s_1}\textrm{Im }F_{l_{k}}...\textrm{Im }F_{l_{1}}\textrm{Re }F_j\textrm{Im }F_{l_1}...\textrm{Im }F_{l_{k}}(\textrm{Im }F_{p})^{s_2} \\
+\frac{1}{2}(-1)^{k+s_2} (\textrm{Im }F_{p})^{s_2}\textrm{Im }F_{l_{k}}...\textrm{Im }F_{l_{1}}\textrm{Re }F_j\textrm{Im }F_{l_1}...\textrm{Im }F_{l_{k}} (\textrm{Im }F_{p})^{s_1},\num
\end{multline*} 
since 
\inc\begin{align*}\label{mari6}
& \ (-1)^{k+s_1}\sigma\big(X,(\textrm{Im }F_{p})^{s_1}\textrm{Im }F_{l_{k}}...\textrm{Im }F_{l_{1}}\textrm{Re }F_j\textrm{Im }F_{l_1}...\textrm{Im }F_{l_{k}}(\textrm{Im }F_{p})^{s_2}X\big)\num\\
= & \ \sigma\big(\textrm{Im }F_{l_1}...\textrm{Im }F_{l_{k}}(\textrm{Im }F_{p})^{s_1}X,\textrm{Re }F_j\textrm{Im }F_{l_1}...\textrm{Im }F_{l_{k}}(\textrm{Im }F_{p})^{s_2}X\big)  \\
= &\ \textrm{Re }q_j\big(\textrm{Im }F_{l_1}...\textrm{Im }F_{l_{k}}(\textrm{Im }F_{p})^{s_1}X;\textrm{Im }F_{l_1}...\textrm{Im }F_{l_{k}}(\textrm{Im }F_{p})^{s_2}X\big) \\
= &\  \textrm{Re }q_j\big(\textrm{Im }F_{l_1}...\textrm{Im }F_{l_{k}}(\textrm{Im }F_{p})^{s_2}X;\textrm{Im }F_{l_1}...\textrm{Im }F_{l_{k}}(\textrm{Im }F_{p})^{s_1}X\big)\\
= & \  \sigma\big(\textrm{Im }F_{l_1}...\textrm{Im }F_{l_{k}}(\textrm{Im }F_{p})^{s_2}X,\textrm{Re }F_j\textrm{Im }F_{l_1}...\textrm{Im }F_{l_{k}}(\textrm{Im }F_{p})^{s_1}X\big)  \\
= & \ (-1)^{k+s_2}\sigma\big(X,(\textrm{Im }F_{p})^{s_2}\textrm{Im }F_{l_{k}}...\textrm{Im }F_{l_{1}}\textrm{Re }F_j\textrm{Im }F_{l_1}...\textrm{Im }F_{l_{k}}(\textrm{Im }F_{p})^{s_1}X\big).
\end{align*}
Then, a direct computation (see Lemma~2 in~\cite{mz}) shows that the Hamilton map of the quadratic form 
$$H_{\textrm{Im}q_p} \tilde{r}=\big\{\textrm{Im }q_p,\tilde{r}\big\}=\frac{\partial \textrm{Im }q_p}{\partial \xi}.\frac{\partial \tilde{r}}{\partial x}-\frac{\partial \textrm{Im } q_p}{\partial x}.\frac{\partial \tilde{r}}{\partial \xi},$$
is given by the commutator $-2[\textrm{Im }F_p,\tilde{F}]$, that is,
$$H_{\textrm{Im}q_p} \tilde{r}(X)=-2\sigma\big(X,[\textrm{Im }F_p,\tilde{F}]X\big).$$
A computation as in (\ref{mari6}) then allows to directly get (\ref{mari5}).~$\Box$

\bigskip

\begin{lemma}\label{lem3}
Consider a $C^{\infty}(\rr)$ function $f$ such that 
$$f \in L^{\infty}(\rr) \textrm{ and } \exists c_1,c_2>0, \ \emph{\textrm{supp }}f' \subset \big\{x \in \rr : c_1 \leq |x| \leq c_2\big\},$$
and $r$ a non-negative quadratic form then for all $0<\alpha \leq 1$,
\begin{equation}\label{giu16}\inc
f\big(r(X)\langle X\rangle^{-2\alpha}\big) \in S(1,\langle X\rangle^{-2\alpha}dX^2). \num
\end{equation}
\end{lemma}

\bigskip

\noindent
\textit{Proof of Lemma~\ref{lem3}}. It is sufficient to check that 
\begin{equation}\label{giu17}\inc
\nabla\big(r(X)\langle X\rangle^{-2\alpha}\big) \in S_{\Omega}\big(\langle X \rangle^{-\alpha},\langle X \rangle^{-2\alpha}dX^2\big), \num
\end{equation}
where $\Omega$ is a small open neighborhood of  $\textrm{supp }f'\big(r(X)\langle X\rangle^{-2\alpha}\big).$ We deduce from (\ref{giu00.1}) and the fact that $r(X)$ is a non-negative quadratic form that 
$$r(X) \sim \langle X \rangle^{2\alpha}$$
and
$$|\nabla r(X)\big)| \lesssim r(X)^{1/2} \lesssim \langle X\rangle^{\alpha},$$
on $\Omega$. By noticing that $0<\alpha \leq 1$, $\langle X \rangle^{r} \in S(\langle X \rangle^r,\langle X \rangle^{-2}dX^2)$, for any $r \in \rr$; and that the function 
$r(X)$ is just a quadratic form, we directly deduce (\ref{giu17}) from the previous estimates and the Leibniz's rule, since
$$r(X) \in S_{\Omega}\big(\langle X \rangle^{2\alpha},\langle X \rangle^{-2\alpha}dX^2\big). \ \Box$$

\bigskip

In all the following lemmas, we shall denote by $r_k$ the quadratic forms defined in (\ref{bellagiu11}) for $0 \leq k \leq m$.

\bigskip

\begin{lemma}\label{lem2.15}
For all $s \in \rr$ and $0 \leq j \leq m-2$, we have
$$r_{m-j-1}(X)^{s} \in S_{\Omega}\big(r_{m-j-1}(X)^{s},r_{m-j-1}(X)^{-1}dX^2\big),$$
if $\Omega$ is any open set where
$$r_{m-j-1}(X) \gtrsim \langle X \rangle^{\frac{2(2m-2j-1)}{2m+1}}.$$  
\end{lemma}

\bigskip

\noindent
\textit{Proof of Lemma~\ref{lem2.15}}. Recalling from (\ref{bellagiu11}) that the symbol $r_{m-j-1}(X)$ is a non-negative quadratic form and that we have from (\ref{giu00.1}) that 
\begin{equation}\label{mar1}\inc
|\nabla r_{m-j-1}(X)| \lesssim r_{m-j-1}(X)^{\frac{1}{2}}, \num
\end{equation}
which implies that for all $s \in \rr$,
\inc\begin{align*}\label{na3}
& \ \frac{\big|\nabla\big(r_{m-j-1}(X)^{s}\big)\big|}{r_{m-j-1}(X)^{s}}
\lesssim \frac{\big|\nabla r_{m-j-1}(X)\big|}{r_{m-j-1}(X)} \num \\
\lesssim & \ r_{m-j-1}(X)^{-\frac{1}{2}},
\end{align*}
on $\Omega$, we notice that the result of Lemma~\ref{lem2.15} is therefore a straightforward consequence of the Leibniz's rule.~$\Box$

\bigskip

\begin{lemma}\label{lem2.1}
Consider the function $\Psi_j$ defined in \emph{(\ref{giu39})} then for any $0 \leq j \leq m-2$,
$$\Psi_j \in S_{\Omega}\Big(1,r_{m-j-1}(X)^{-\frac{2m-2j-3}{2m-2j-1}}dX^2\Big),$$
if $\Omega$ is any open set where
$$r_{m-j-1}(X) \gtrsim \langle X \rangle^{\frac{2(2m-2j-1)}{2m+1}},$$  
which implies in particular that
$$\Psi_j \in S_{\Omega}\big(1,\langle X \rangle^{-\frac{2(2m-2j-3)}{2m+1}}dX^2\big).$$
\end{lemma}
\bigskip

\noindent
\textit{Proof of Lemma~\ref{lem2.1}}. We first notice from (\ref{giu13}) and (\ref{giu39}) that
$$r_{m-j-2}(X) \sim r_{m-j-1}(X)^{\frac{2m-2j-3}{2m-2j-1}},$$
on $\Omega \cap \textrm{supp }\Psi_j'$. 
Since from (\ref{giu00.1}), 
\inc\begin{align*}\label{na1}
|\nabla r_{m-j-2}(X)| \lesssim & \ r_{m-j-2}(X)^{\frac{1}{2}} \num \\
\lesssim & \  r_{m-j-1}(X)^{\frac{2m-2j-3}{2(2m-2j-1)}},
\end{align*}
on $\Omega \cap \textrm{supp }\Psi_j'$, we deduce that the quadratic symbol $r_{m-j-2}(X)$ belongs to the class
\begin{equation}\label{na2}\inc
S_{\Omega \cap \textrm{supp}\Psi_j'}\Big(r_{m-j-1}(X)^{\frac{2m-2j-3}{2m-2j-1}},\frac{dX^2}{r_{m-j-1}(X)^{\frac{2m-2j-3}{2m-2j-1}}}\Big).\num
\end{equation}
It follows from Lemma~\ref{lem2.15} that 
$$
\frac{r_{m-j-2}(X)}{r_{m-j-1}(X)^{\frac{2m-2j-3}{2m-2j-1}}}
\in S_{\Omega \cap \textrm{supp}\Psi_j'}\Big(1,\frac{dX^2}{r_{m-j-1}(X)^{\frac{2m-2j-3}{2m-2j-1}}}\Big),$$
which implies that 
$$\Psi_j \in S_{\Omega}\big(1,r_{m-j-1}(X)^{-\frac{2m-2j-3}{2m-2j-1}}dX^2\big).$$
This ends the proof of Lemma~\ref{lem2.1}.~$\Box$

\bigskip

\begin{lemma}\label{lem2.2}
Consider the function $W_j$ defined in \emph{(\ref{giu40})} then for any $1 \leq j \leq m-1$,
$$W_j \in S_{\Omega}\big(1,r_{m-j-1}(X)^{-1}dX^2\big),$$
if $\Omega$ is any open set where
$$r_{m-j-1}(X) \gtrsim \langle X \rangle^{\frac{2(2m-2j-1)}{2m+1}},$$  
which implies in particular that 
$$W_j \in S_{\Omega}\big(1,\langle X \rangle^{-\frac{2(2m-2j-1)}{2m+1}}dX^2\big).$$
\end{lemma}

\bigskip

\noindent
\textit{Proof of Lemma~\ref{lem2.2}}. By noticing from (\ref{giu15}) and (\ref{giu40}) that
$$r_{m-j-1}(X) \sim r_{m-j}(X)^{\frac{2m-2j-1}{2m-2j+1}}$$
and
$$r_{m-j}(X) \gtrsim \langle X \rangle^{\frac{2(2m-2j+1)}{2m+1}},$$
on $\Omega \cap \textrm{supp }W_j'$, and that the two derivatives $\psi'$ and $w_2'$ of the functions appearing in (\ref{giu39}) and (\ref{giu40}) have similar types of support as the function defined in (\ref{giu14}), we notice that we are exactly in the setting studied in Lemma~\ref{lem2.1} with $j$ replaced by $j-1$. We therefore deduce the result of Lemma~\ref{lem2.2} from our analysis led in the proof of Lemma~\ref{lem2.1}.~$\Box$

\bigskip

\begin{lemma}\label{lem2.21}
If $s_1$, $s_2 \in \nn$, $1 \leq j,p \leq N$, $(l_1,...,l_k) \in \{1,...,N\}^k$ then we have 
\begin{align*}
& \ \big|\emph{\textrm{Re }}q_j(\emph{\textrm{Im }}F_{l_1}...\emph{\textrm{Im }}F_{l_{k}}(\emph{\textrm{Im }}F_{p})^{s_1}X;\emph{\textrm{Im }}F_{l_1}...\emph{\textrm{Im }}F_{l_{k}}(\emph{\textrm{Im }}F_{p})^{s_2}X)\big| \\
\leq & \ \emph{\textrm{Re }}q_j(\emph{\textrm{Im }}F_{l_1}...\emph{\textrm{Im }}F_{l_{k}}(\emph{\textrm{Im }}F_{p})^{s_1}X)^{\frac{1}{2}}
\emph{\textrm{Re }}q_j(\emph{\textrm{Im }}F_{l_1}...\emph{\textrm{Im }}F_{l_{k}}(\emph{\textrm{Im }}F_{p})^{s_2}X)^{\frac{1}{2}}\\
\leq & \ r_{k+s_1}(X)^{\frac{1}{2}}r_{k+s_2}(X)^{\frac{1}{2}}
\end{align*}
and
\begin{align*}
& \ \big|\nabla \big[ \emph{\textrm{Re }}q_j(\emph{\textrm{Im }}F_{l_1}...\emph{\textrm{Im }}F_{l_{k}}(\emph{\textrm{Im }}F_{p})^{s_1}X;\emph{\textrm{Im }}F_{l_1}...\emph{\textrm{Im }}F_{l_{k}}(\emph{\textrm{Im }}F_{p})^{s_2}X)\big]\big| \\
\lesssim & \ \emph{\textrm{Re }}q_j(\emph{\textrm{Im }}F_{l_1}...\emph{\textrm{Im }}F_{l_{k}}(\emph{\textrm{Im }}F_{p})^{s_1}X)^{\frac{1}{2}}+\emph{\textrm{Re }}q_j(\emph{\textrm{Im }}F_{l_1}...\emph{\textrm{Im }}F_{l_{k}}(\emph{\textrm{Im }}F_{p})^{s_2}X)^{\frac{1}{2}} \\
\lesssim & \ r_{k+\emph{\textrm{max}}(s_1,s_2)}(X)^{\frac{1}{2}}.
\end{align*}
\end{lemma}

\bigskip

\noindent
\textit{Proof of Lemma~\ref{lem2.21}}. By reason of symmetry, we can assume in the following that $s_1 \leq s_2$. Recalling that the quadratic form $\textrm{Re }q_j$ is non-negative, the first estimate is a direct consequence of (\ref{bellagiu11}) and the Cauchy-Schwarz inequality.
About the second estimate, we recall from (\ref{lol10}) that the Hamilton map of the quadratic form
$$\textrm{Re }q_j(\textrm{Im }F_{l_1}...\textrm{Im }F_{l_{k}}(\textrm{Im }F_{p})^{s_1}X;\textrm{Im }F_{l_1}...\textrm{Im }F_{l_{k}}(\textrm{Im }F_{p})^{s_2}X),$$ 
is 
\begin{multline*}
\frac{1}{2}(-1)^{k+s_1}(\textrm{Im }F_{p})^{s_1}\textrm{Im }F_{l_{k}}...\textrm{Im }F_{l_{1}}\textrm{Re }F_j\textrm{Im }F_{l_1}...\textrm{Im }F_{l_{k}}(\textrm{Im }F_{p})^{s_2} \\
+\frac{1}{2}(-1)^{k+s_2} (\textrm{Im }F_{p})^{s_2}\textrm{Im }F_{l_{k}}...\textrm{Im }F_{l_{1}}\textrm{Re }F_j\textrm{Im }F_{l_1}...\textrm{Im }F_{l_{k}} (\textrm{Im }F_{p})^{s_1}.
\end{multline*} 
A direct computation as in (3.18) of~\cite{mz} shows that
\inc\begin{align*}\label{mari1}
& \ \num \nabla \big[\textrm{Re }q_j(\textrm{Im }F_{l_1}...\textrm{Im }F_{l_{k}}(\textrm{Im }F_{p})^{s_1}X;\textrm{Im }F_{l_1}...\textrm{Im }F_{l_{k}}(\textrm{Im }F_{p})^{s_2}X)\big]\\
= & \ (-1)^{k+s_1+1}\sigma(\textrm{Im }F_{p})^{s_1}\textrm{Im }F_{l_{k}}...\textrm{Im }F_{l_{1}}\textrm{Re }F_j\textrm{Im }F_{l_1}...\textrm{Im }F_{l_{k}}(\textrm{Im }F_{p})^{s_2} \\
+ & \ (-1)^{k+s_2+1} \sigma(\textrm{Im }F_{p})^{s_2}\textrm{Im }F_{l_{k}}...\textrm{Im }F_{l_{1}}\textrm{Re }F_j\textrm{Im }F_{l_1}...\textrm{Im }F_{l_{k}} (\textrm{Im }F_{p})^{s_1}
\end{align*}
where
$$\sigma=\begin{pmatrix}
0 & I_n\\
-I_n & 0
\end{pmatrix}.
$$ 
The notation $I_n$ stands here for the $n$ by $n$ identity matrix. We deduce from (\ref{giu00.1}) and (\ref{mari1}) that for any $s \in \nn$,
\inc\begin{align*}\label{kee4}
& \ |(\textrm{Im }F_{p})^{s}\textrm{Im }F_{l_{k}}...\textrm{Im }F_{l_{1}}\textrm{Re }F_j\textrm{Im }F_{l_1}...\textrm{Im }F_{l_{k}}(\textrm{Im }F_{p})^{s}X|\num \\
\lesssim & \  \big|\nabla \big[\textrm{Re }q_j(\textrm{Im }F_{l_1}...\textrm{Im }F_{l_{k}}(\textrm{Im }F_{p})^{s}X)\big] \big| \\
\lesssim & \  \textrm{Re }q_j(\textrm{Im }F_{l_1}...\textrm{Im }F_{l_{k}}(\textrm{Im }F_{p})^{s}X)^{\frac{1}{2}}.
\end{align*}
By using twice the estimate (\ref{kee4}) with respectively $X$ and $(\textrm{Im }F_p)^{s_2-s_1}X$, and the index $s=s_1$, we deduce from (\ref{bellagiu11}) and (\ref{mari1}) the second estimate in Lemma~\ref{lem2.21}.~$\Box$

\bigskip

\begin{lemma}\label{lem2.3}
Consider the quadratic form $\tilde{r}_{m-j-1,p}$ defined in \emph{(\ref{giu11})} and \emph{(\ref{giu11serena})} then for any $0 \leq j \leq m-2$ and $1 \leq p \leq N$,
$$\frac{\tilde{r}_{m-j-1,p}(X)}{r_{m-j-1}(X)^{\frac{2m-2j-2}{2m-2j-1}}} \in 
S_{\Omega}\big(1,r_{m-j-1}(X)^{-\frac{2m-2j-3}{2m-2j-1}}dX^2\big),$$
if $\Omega$ is any open set where
$$r_{m-j-1}(X) \gtrsim \langle X \rangle^{\frac{2(2m-2j-1)}{2m+1}}$$ 
and  
$$r_{m-j-2}(X) \lesssim r_{m-j-1}(X)^{\frac{2m-2j-3}{2m-2j-1}},$$  
which implies in particular that 
$$\frac{\tilde{r}_{m-j-1,p}(X)}{r_{m-j-1}(X)^{\frac{2m-2j-2}{2m-2j-1}}} \in S_{\Omega}\big(1,\langle X \rangle^{-\frac{2(2m-2j-3)}{2m+1}}dX^2\big).$$
\end{lemma}

\bigskip

\noindent
\textit{Proof of Lemma~\ref{lem2.3}}. Since from Lemma~\ref{lem2.21},
$$|\tilde{r}_{m-j-1,p}(X)| \lesssim r_{m-j-1}(X)^{\frac{2m-2j-2}{2m-2j-1}}$$
and
\begin{align*}
|\nabla \tilde{r}_{m-j-1,p}(X)| \lesssim & \ r_{m-j-1}(X)^{\frac{1}{2}}+r_{m-j-2}(X)^{\frac{1}{2}} \\
\lesssim & \  r_{m-j-1}(X)^{\frac{1}{2}},
\end{align*} 
on $\Omega$, we get that the quadratic form $\tilde{r}_{m-j-1,p} $ belongs to the symbol class
$$S_{\Omega}\big(r_{m-j-1}(X)^{\frac{2m-2j-2}{2m-2j-1}},r_{m-j-1}(X)^{-\frac{2m-2j-3}{2m-2j-1}}dX^2\big).$$
One can then deduce the result of Lemma~\ref{lem2.3} from Lemma~\ref{lem2.15}.~$\Box$

\bigskip

When adding a large parameter $\Lambda_j \geq 1$ in the description of the open set $\Omega$, a straightforward adaptation of the proof of the previous lemma gives the 
following $L^{\infty}(\Omega)$ estimate with respect to this parameter.

\bigskip

\begin{lemma}\label{lem2.331}
Consider the quadratic form $\tilde{r}_{m-j-1,p}$ defined in \emph{(\ref{giu11})} and \emph{(\ref{giu11serena})} then for any $0 \leq j \leq m-2$ and $1 \leq p \leq N$,
$$\big\|r_{m-j-1}(X)^{-\frac{2m-2j-2}{2m-2j-1}}\tilde{r}_{m-j-1,p}(X)   \big\|_{L^{\infty}(\Omega)} \lesssim \Lambda_j^{-\frac{1}{2}},$$
if $\Omega$ is any open set where
$$r_{m-j-1}(X) \gtrsim \langle X \rangle^{\frac{2(2m-2j-1)}{2m+1}}$$ 
and  
$$r_{m-j-2}(X) \lesssim \Lambda_j^{-1} r_{m-j-1}(X)^{\frac{2m-2j-3}{2m-2j-1}},$$
with  $\Lambda_j \geq 1$.  
\end{lemma}

\bigskip

In the following lemmas, we shall carefully study the dependence of the estimates with respect to the large parameter $\Lambda_j \geq 1$.

\bigskip

\begin{lemma}\label{lem2.4}
For any $0 \leq j \leq m-2$, we have for all $X \in \Omega$,
$$\left|\sum_{p=1}^N\frac{H_{\emph{\textrm{Im}}q_p}\tilde{r}_{m-j-1,p}(X)}{r_{m-j-1}(X)^{\frac{2m-2j-2}{2m-2j-1}}} 
-2r_{m-j-1}(X)^{\frac{1}{2m-2j-1}}\right| \\
\lesssim \Lambda_j^{-\frac{1}{2}} r_{m-j-1}(X)^{\frac{1}{2m-2j-1}},$$
if $\Omega$ is any open set where
$$r_{m-j-1}(X) \gtrsim \langle X \rangle^{\frac{2(2m-2j-1)}{2m+1}},$$
$$r_{m-j-2}(X) \lesssim \Lambda_j^{-1} r_{m-j-1}(X)^{\frac{2m-2j-3}{2m-2j-1}},$$
$$r_{m-j}(X) \lesssim r_{m-j-1}(X)^{\frac{2m-2j+1}{2m-2j-1}},$$
with  $\Lambda_j \geq 1$.
\end{lemma}
\bigskip

\noindent
\textit{Proof of Lemma~\ref{lem2.4}}. We begin by writing from (\ref{giu11}), (\ref{giu11serena}) and Lemma~\ref{lem2} that 
\inc\begin{multline*} \label{kee1}
H_{\textrm{Im}q_p}\tilde{r}_{m-j-1,p}(X)=
2\hspace{-0.5cm} \sum_{\substack{s=1,...,N \\ (l_1,...,l_{m-j-2}) \in \{1,...,N\}^{m-j-2}}}\hspace{-0.5cm} \textrm{Re }q_s(\textrm{Im }F_{l_1}...\textrm{Im }F_{l_{m-j-2}}\textrm{Im }F_{p}X)
\\ +  2\hspace{-0.5cm} \sum_{\substack{s=1,...,N \\ (l_1,...,l_{m-j-2}) \in \{1,...,N\}^{m-j-2}}} \hspace{-0.5cm} \textrm{Re }q_s(\textrm{Im }F_{l_1}...\textrm{Im }F_{l_{m-j-2}}X;\textrm{Im }F_{l_1}...\textrm{Im }F_{l_{m-j-2}}(\textrm{Im }F_{p})^2X).\num
\end{multline*}
Lemma~\ref{lem2.4} is then a consequence of the following estimate 
\begin{align*}
& \ \big|    \textrm{Re }q_s(\textrm{Im }F_{l_1}...\textrm{Im }F_{l_{m-j-2}}X;\textrm{Im }F_{l_1}...\textrm{Im }F_{l_{m-j-2}}(\textrm{Im }F_{p})^2X)\big| \\
\leq & \ \textrm{Re }q_s(\textrm{Im }F_{l_1}...\textrm{Im }F_{l_{m-j-2}}X)^{\frac{1}{2}}
\textrm{Re }q_s(\textrm{Im }F_{l_1}...\textrm{Im }F_{l_{m-j-2}}(\textrm{Im }F_{p})^2X)^{\frac{1}{2}}\\
\leq & \ r_{m-j-2}(X)^{\frac{1}{2}}r_{m-j}(X)^{\frac{1}{2}}\\
\lesssim & \ \Lambda_j^{-\frac{1}{2}}r_{m-j-1}(X),
\end{align*}
fulfilled on $\Omega$ that we obtain from the Cauchy-Schwarz inequality.~$\Box$

\bigskip

\begin{lemma}\label{lem2.41}
For any $0 \leq j \leq m-2$ and $1 \leq p \leq N$, we have for all $X \in \Omega$,
$$\big|r_{m-j-1}(X)^{\frac{2m-2j-2}{2m-2j-1}}H_{\emph{\textrm{Im}}q_p}\big(r_{m-j-1}(X)^{-\frac{2m-2j-2}{2m-2j-1}} \big)\big|  
\lesssim r_{m-j-1}(X)^{\frac{1}{2m-2j-1}},$$
if $\Omega$ is any open set where
$$r_{m-j-1}(X) \gtrsim \langle X \rangle^{\frac{2(2m-2j-1)}{2m+1}},$$
$$r_{m-j-2}(X) \lesssim \Lambda_j^{-1}r_{m-j-1}(X)^{\frac{2m-2j-3}{2m-2j-1}},$$
$$r_{m-j}(X) \lesssim r_{m-j-1}(X)^{\frac{2m-2j+1}{2m-2j-1}},$$
with  $\Lambda_j \geq 1$.
\end{lemma}
\bigskip

\noindent
\textit{Proof of Lemma~\ref{lem2.41}}. We begin by writing from (\ref{bellagiu11}) and Lemma~\ref{lem2} that 
\inc\begin{multline*} \label{kee1.001}
H_{\textrm{Im}q_p} r_{m-j-1}(X) \\
=4\hspace{-0.5cm}\sum_{\substack{s=1,...,N \\ (l_1,...,l_{m-j-1}) \in \{1,...,N\}^{m-j-1}}}\hspace{-0.5cm} \textrm{Re }q_s(\textrm{Im }F_{l_1}...\textrm{Im }F_{l_{m-j-1}}X;\textrm{Im }F_{l_1}...\textrm{Im }F_{l_{m-j-1}}\textrm{Im }F_p X). \num
\end{multline*}
Since 
\begin{align*}
& \ r_{m-j-1}(X)^{\frac{2m-2j-2}{2m-2j-1}}H_{\textrm{Im}q_p}\Big(r_{m-j-1}(X)^{-\frac{2m-2j-2}{2m-2j-1}} \Big)\\
= & \ -\frac{2m-2j-2}{2m-2j-1}\frac{H_{\textrm{Im}q_p}r_{m-j-1}(X)}{r_{m-j-1}(X)},
\end{align*}
Lemma~\ref{lem2.41} is then a consequence of the following estimate
\inc\begin{align*}\label{spea2}
& \ \num \big| \textrm{Re }q_s(\textrm{Im }F_{l_1}...\textrm{Im }F_{l_{m-j-1}}X;\textrm{Im }F_{l_1}...\textrm{Im }F_{l_{m-j-1}}\textrm{Im }F_p X)  \big| \\
\leq & \ \textrm{Re }q_s(\textrm{Im }F_{l_1}...\textrm{Im }F_{l_{m-j-1}}X)^{\frac{1}{2}}
\textrm{Re }q_s(\textrm{Im }F_{l_1}...\textrm{Im }F_{l_{m-j-1}}\textrm{Im }F_p X)^{\frac{1}{2}}\\
\leq & \ r_{m-j-1}(X)^{\frac{1}{2}} r_{m-j}(X)^{\frac{1}{2}}\\
\lesssim & \ r_{m-j-1}(X)^{1+\frac{1}{2m-2j-1}}, 
\end{align*}
fulfilled on $\Omega$ that we obtain from the Cauchy-Schwarz inequality.~$\Box$

\bigskip

\begin{lemma}\label{lem2.42}
Consider the functions $\Psi_j$ and $W_{j+1}$ defined in \emph{(\ref{giu39})} and \emph{(\ref{giu40})} then for any $0 \leq j \leq m-2$ and $1 \leq p \leq N$, we have for all $X \in \Omega$,
$$|H_{\emph{\textrm{Im}}q_p}\Psi_j(X)| \lesssim \Lambda_j^{\frac{1}{2}} r_{m-j-1}(X)^{\frac{1}{2m-2j-1}}W_{j+1}(X),$$
if $\Omega$ is any open set where
$$r_{m-j-1}(X) \gtrsim \langle X \rangle^{\frac{2(2m-2j-1)}{2m+1}},$$
$$r_{m-j-2}(X) \lesssim \Lambda_j^{-1} r_{m-j-1}(X)^{\frac{2m-2j-3}{2m-2j-1}},$$
$$r_{m-j}(X) \lesssim r_{m-j-1}(X)^{\frac{2m-2j+1}{2m-2j-1}},$$
with  $\Lambda_j \geq 1$.
\end{lemma}
\bigskip

\noindent
\textit{Proof of Lemma~\ref{lem2.42}}. We begin by noticing from (\ref{giu40}) and (\ref{bell1}) that
\begin{equation}\label{spea10}\inc
\left|\psi'\left(\frac{\Lambda_j r_{m-j-2}(X)}{r_{m-j-1}(X)^{\frac{2m-2j-3}{2m-2j-1}}}\right)\right| \lesssim
W_{j+1}(X),\num
\end{equation}
and by writing from Lemma~\ref{lem2} that 
\inc\begin{multline*} \label{kee1.002}
H_{\textrm{Im}q_p}r_{m-j-2}(X) \\ =4\hspace{-0.5cm}\sum_{\substack{s=1,...,N \\ (l_1,...,l_{m-j-2}) \in \{1,...,N\}^{m-j-2}}}\hspace{-0.5cm} \textrm{Re }q_s(\textrm{Im }F_{l_1}...\textrm{Im }F_{l_{m-j-2}}X;\textrm{Im }F_{l_1}...\textrm{Im }F_{l_{m-j-2}}\textrm{Im }F_p X). \num
\end{multline*}
It follows from the Cauchy-Schwarz inequality that for all $X \in \Omega$,
\inc\begin{align*}\label{spea1}
& \ \num \big|\textrm{Re }q_s(\textrm{Im }F_{l_1}...\textrm{Im }F_{l_{m-j-2}}X;\textrm{Im }F_{l_1}...\textrm{Im }F_{l_{m-j-2}}\textrm{Im }F_p X)\big| \\
\leq & \ \textrm{Re }q_s(\textrm{Im }F_{l_1}...\textrm{Im }F_{l_{m-j-2}}X)^{\frac{1}{2}}
\textrm{Re }q_s(\textrm{Im }F_{l_1}...\textrm{Im }F_{l_{m-j-2}}\textrm{Im }F_p X)^{\frac{1}{2}}\\
\leq & \ r_{m-j-2}(X)^{\frac{1}{2}} r_{m-j-1} (X)^{\frac{1}{2}}\\
\lesssim & \ \Lambda_j^{-\frac{1}{2}}r_{m-j-1}(X)^{\frac{2m-2j-2}{2m-2j-1}}. 
\end{align*}
Then, by writing that 
\begin{multline*}
H_{\textrm{Im}q_p}\left(\frac{\Lambda_j r_{m-j-2}(X)}{r_{m-j-1}(X)^{\frac{2m-2j-3}{2m-2j-1}}}\right)=
\frac{\Lambda_j H_{\textrm{Im}q_p}r_{m-j-2}(X)}{r_{m-j-1}(X)^{\frac{2m-2j-3}{2m-2j-1}}}
\\ -\frac{2m-2j-3}{2m-2j-1}\frac{\Lambda_j r_{m-j-2}(X)H_{\textrm{Im}q_p}r_{m-j-1}(X)}{r_{m-j-1}(X)^{1+\frac{2m-2j-3}{2m-2j-1}}}.
\end{multline*}
Lemma~\ref{lem2.42} is a consequence of (\ref{giu39}), (\ref{kee1.001}), (\ref{spea2}), (\ref{kee1.002}), (\ref{spea1}) and (\ref{spea10}), since 
$$r_{m-j-2}(X) \sim \Lambda_j^{-1} r_{m-j-1}(X)^{\frac{2m-2j-3}{2m-2j-1}},$$
on the support of $\Psi_j'$.~$\Box$

\bigskip

\begin{lemma}\label{lem2.435}
For $m \geq 2$, consider the function $\tilde{W}_0$ defined in \emph{(\ref{giu41})} then for all $X \in \rr^{2n}$ and $1 \leq p \leq N$,
$$|H_{\emph{\textrm{Im}}q_p}  \tilde{W}_0(X)| \lesssim \langle X \rangle^{\frac{2}{2m+1}}.$$
\end{lemma}
\bigskip

\noindent
\textit{Proof of Lemma~\ref{lem2.435}}. 
Since $|\nabla \textrm{Im }q_p(X)| \lesssim \langle X \rangle$, because $\textrm{Im }q_p$ is a quadratic form, Lemma~\ref{lem2.435} is then a consequence of (\ref{giu15}), (\ref{bellagiu11}), (\ref{giu41}) and Lemma~\ref{lem3}.~$\Box$

\bigskip

\begin{lemma}\label{lem2.43}
Consider the function $W_{j+1}$ defined in \emph{(\ref{giu40})} then for any $0 \leq j \leq m-2$ and $1 \leq p \leq N$, we have for all $X \in \Omega$, 
$$|H_{\emph{\textrm{Im}}q_p}W_{j+1}(X)| \lesssim \Lambda_j^{\frac{1}{2}} r_{m-j-1}(X)^{\frac{1}{2m-2j-1}}\Psi_j(X),$$
if $\Omega$ is any open set where
$$r_{m-j-1}(X) \gtrsim \langle X \rangle^{\frac{2(2m-2j-1)}{2m+1}},$$
$$r_{m-j-2}(X) \lesssim \Lambda_j^{-1} r_{m-j-1}(X)^{\frac{2m-2j-3}{2m-2j-1}},$$
$$r_{m-j}(X) \lesssim r_{m-j-1}(X)^{\frac{2m-2j+1}{2m-2j-1}},$$
with  $\Lambda_j \geq 1$.
\end{lemma}
\bigskip

\noindent
\textit{Proof of Lemma~\ref{lem2.43}}. 
One can notice from (\ref{giu13}), (\ref{giu15}), (\ref{giu39}), (\ref{giu40}) and (\ref{bell1}) that
\begin{equation}\label{spea10}\inc
\forall \ 0 \leq j \leq m-2, \ \left|w_2'\left(\frac{\Lambda_j r_{m-j-2}(X)}{r_{m-j-1}(X)^{\frac{2m-2j-3}{2m-2j-1}}}\right)\right| \lesssim
\Psi_{j}(X),\num
\end{equation}
and that the derivatives of  $\Psi_j$ and $W_{j+1}$ are exactly the same types of functions. It follows that Lemma~\ref{lem2.43} is just a straightforward consequence of Lemma~\ref{lem2.42}.~$\Box$

\end{document}